\documentclass[11pt,reqno,centertags]{amsart}
\usepackage{amsmath,latexsym}
\usepackage[psamsfonts]{amssymb}
\usepackage{amssymb}
\usepackage{upref}
\usepackage{amsmath}
\usepackage{mathptmx}
\usepackage{mathrsfs}

\DeclareSymbolFont{SY}{U}{psy}{m}{n}
\DeclareMathSymbol{\emptyset}{\mathord}{SY}{'306}

\newcommand{\bbC}{\mathbb{C}}
\newcommand{\bbN}{\mathbb{N}}
\newcommand{\bbR}{\mathbb{R}}
\newcommand{\cD}{{\mathcal D}}

\newcommand{\cW}{{\mathcal W}}
\newcommand{\sE}{\mathsf{E}}
\newcommand{\fH}{\mathfrak{H}}

\newcommand{\fM}{\mathfrak{M}}
\newcommand{\fN}{\mathfrak{N}}

\newcommand{\lal}{\langle}
\newcommand{\ral}{\rangle}
\newcommand{\got}{\mathfrak}

\newcommand{\im}{\mathop{\mathrm{Im}}}
\newcommand{\re}{\mathop{\mathrm{Re}}}
\newcommand{\ran}{\mathop{\mathrm{Ran}}}
\newcommand{\Ran}{\mathop{\mathrm{Ran}}}
\renewcommand{\ker}{{\ensuremath{\mathrm{Ker}}}}
\newcommand{\Ker}{{\ensuremath{\mathrm{Ker}}}}

\newcommand{\spec}{\mathop{\text{\rm spec}}}
\newcommand{\dist}{\mathop{\mathrm{dist}}}
\newcommand{\dom}{\mathop{\mathrm{Dom}}}
\newcommand{\Dom}{\mathop{\mathrm{Dom}}}

\newtheorem{introtheorem}{Theorem}{\bf}{\it}

\newtheorem{thm}{Theorem}[section]
\newtheorem{hypo}[thm]{Hypothesis}
\newtheorem{cor}[thm]{Corollary}
\newtheorem{lem}[thm]{Lemma}

\theoremstyle{definition}
\newtheorem{defn}[thm]{Definition}
\theoremstyle{remark}
\newtheorem{rem}[thm]{Remark}
\newtheorem{ex}[thm]{Example}
\newtheorem{introremark}[introtheorem]{Remark}{\it}{\rm}

\numberwithin{equation}{section}

\begin{document}

\title[\qquad Sharp norm estimates in the subspace perturbation
problem] {Some sharp norm estimates in the subspace \\ perturbation
problem}

\author[A. K. Motovilov]{Alexander K. Motovilov}

\address{A. K. Motovilov\\
Bogoliubov Laboratory of Theoretical Physics\\
Joint Institute for Nuclear Research\\
141980 Dubna, Moscow Region, Russia}

\email{motovilv@theor.jinr.ru}

\author[A. V. Selin]{Alexei V. Selin}
\address{A. V. Selin\\ Laboratory of Informational Technologies\\
Joint Institute for Nuclear Research\\
141980 Dubna, Moscow Region, Russia}

\email{selin@theor.jinr.ru}

\subjclass{Primary 47A55; Secondary 47B25}

\keywords{Perturbation problem, spectral subspaces, direct rotation,
numerical range}

\date{February 27, 2005}

\begin{abstract}
We discuss the spectral subspace perturbation problem for a
self-adjoint operator. Assuming that the convex hull of a part of
its spectrum does not intersect the remainder of the spectrum, we
establish an \textit{a priori} sharp bound on variation of the
corresponding spectral subspace under off-diagonal perturbations.
This bound represents a new, \textit{a priori}, $\tan\Theta$
Theorem. We also extend the Davis--Kahan $\tan 2\Theta$ Theorem to
the case of some unbounded perturbations.
\end{abstract}

\maketitle

\section{Introduction}

Assume that the spectrum of a self-adjoint operator $A$ on a Hilbert
space $\fH$ consists of two disjoint components $\sigma_{-}$ and
$\sigma_{+}$, i.e. $\spec(A) = \sigma_{-} \cup \sigma_{+}$ and
\begin{equation}
\label{ddist} d=\dist(\sigma_{-}, \sigma_{+})  > 0.
\end{equation}
Then $\got H$ is decomposed into the orthogonal sum $\got H = \got
H_{-}\oplus\got H_{+}$ of the spectral subspaces $\got H_{\pm} =
\ran \sE_{A}(\sigma_{\pm})$ where $\sE_A(\delta)$ denotes the
spectral projection of $A$ associated with a Borel set
$\delta\subset\bbR$. It is well known (see, e.g., \cite[\S
135]{SzNagy}) that sufficiently small self-adjoint perturbation $V$
of $A$ does not close the gaps between the sets $\sigma_{-}$ and
$\sigma_{+}$ which allows one to think of the corresponding disjoint
spectral components $\sigma'_{-}$ and $\sigma'_{+}$ of the perturbed
operator $L=A + V$ as a result of the perturbation of the spectral
sets $\sigma_{-}$ and $\sigma_{+}$, respectively. Moreover, the
decomposition $\got H = \got H'_{-}\oplus\got H'_{+}$ with $\got
H'_{\pm} = \ran \sE_{L}(\sigma'_{\pm})$ is continuous in $V$ in the
sense that the projections $\sE_{L}(\sigma'_{\pm})$ converge to
$\sE_{A}(\sigma_{\pm})$ in the operator norm topology as $\|V\| \to
0$.

Given a mutual disposition of the spectral components $\sigma_{\pm}$
of the operator $A$, the problem of perturbation theory is to study
variation of these components and the corresponding spectral
subspaces under the perturbation $V$. In particular, the questions
of interest are as follows (see \cite{KMM1}, \cite{KMM4}):

$(i)$ Under what (sharp) condition on $\|V\|$ do the gaps between
the sets $\sigma_-$ and $\sigma_+$ remain open, i.e.
$\dist(\sigma'_-,\sigma'_+)>0$?

$(ii)$ Having established this condition, can one ensure that it
implies inequality
\begin{equation}
\label{edif} \| \sE_{L}(\sigma'_{-}) - \sE_{A}(\sigma_{-}) \| <1\,?
\end{equation}
(Surely, \eqref{edif} holds if and only if inequality
$\|\sE_{L}(\sigma'_{+}) - \sE_{A}(\sigma_{+})\|<1$ does.)

In general, answer to the question $(i)$ is well known: the gaps
between $\sigma_-$ and $\sigma_+$ remain open if
\begin{equation}
\label{Vdgen} \|V\|<\frac{d}{2}.
\end{equation}
Among all perturbations of the operator $A$ we distinguish the ones
that are off-diagonal with respect to the decomposition $\fH=\ran
\sE_A(\sigma_{-})\oplus\ran \sE_A(\sigma_{+})$, i.e. the
perturbations that anticommute with the difference
\begin{equation}
\label{Jintro} J=\sE_A(\sigma_+)-\sE_A(\sigma_-)
\end{equation}
of the spectral projections $\sE_A(\sigma_+)$ and $\sE_A(\sigma_-)$.
If one restricts oneself to perturbations $V$ of this class then
inequality $\dist(\sigma'_-,\sigma'_+)>0$ is ensured by the weaker
condition
\begin{equation}
\label{sqrt32} \| V \| < \frac{\sqrt{3}}{2} d
\end{equation}
proven in \cite[Theorem 1]{KMM4}. Similarly to \eqref{Vdgen},
condition \eqref{sqrt32} is sharp.

For a review of the known answers to the question $(ii)$ we refer to
\cite{KMM1} in case of the general bounded perturbations and to
\cite{KMM4} in case of the off-diagonal ones. Notice that complete
answers to the question $(ii)$ were found only under certain
additional assumptions on the mutual disposition of the sets
$\sigma_-$ and $\sigma_+$. It is still an open problem whether or
not the corresponding conditions \eqref{Vdgen} and \eqref{sqrt32}
imply \eqref{edif} under the only assumption \eqref{ddist}.

In the present paper we are concerned with the off-diagonal
perturbations and restrict ourselves to two particular mutual
dispositions of the spectral sets $\sigma_-$ and $\sigma_+$. The
first one corresponds to the case where the sets $\sigma_-$ and
$\sigma_+$ are subordinated, say
\begin{equation}
\label{tan2t} \sup \sigma_{-} <  \inf \sigma_{+}.
\end{equation}
The second case under consideration corresponds to a disposition
with one of the sets $\sigma_-$ and $\sigma_+$ lying in a (finite) gap
of the other set, say
\begin{equation}
\label{tant} \sigma_{+} \cap\mathop{\mathrm{conv}}( \sigma_{-} ) =
\emptyset,
\end{equation}
where $\mathop{\mathrm{conv}}(\sigma)$ denotes the convex hull of a
set $\sigma\subset\bbR$.

In both these cases the perturbed spectral sets $\sigma'_{-}$ and
$\sigma'_{+}$ are known to remain disjoint under requirements on
$\|V\|$ much weaker than that of \eqref{sqrt32}.

In particular, if \eqref{tan2t} holds then for any bounded
off-diagonal perturbation $V$ the interval
$(\sup\sigma_-,\inf\sigma_+)$ is in the resolvent set of the
perturbed operator $L=A+V$, and thus $\sigma'_{-} \subset (-\infty,
\sup \sigma_{-}]$ and $\sigma'_{+} \subset [\inf \sigma_{+},
+\infty)$ (see \cite{ALT}, \cite{DK}; cf. \cite{KMM3}). Moreover, in
this case the following norm estimate holds \cite{DK}
\begin{equation*}
%\label{tan2t1}
\| \sE_{L}(\sigma'_{-}) - \sE_{A}(\sigma_{-}) \| \le \sin \Bigl(
\frac{1}{2} \arctan \frac{2\|V\|}{d}\Bigr) < \frac{\sqrt{2}}{2}.
\end{equation*}
This (sharp) bound on the difference of the spectral projection
$\sE_{L}(\sigma'_{-})$ and $\sE_{A}(\sigma_{-})$ is known as the
Davis--Kahan $\tan 2\Theta$ Theorem since it can be written in the
equivalent form $\|\tan 2\Theta \| \le \dfrac{\|V \|}{d}$ where
$\Theta$ is the operator angle between the subspaces $\got H_{-}'$
and $\got H_{-}$ (or between the subspaces $\got H_{+}'$ and $\got
H_{+}$). For definition of the operator angle between two subspaces
see, e.g., \cite{KMM2}.

Our first principal result is an extension of the $\tan 2\Theta$
Theorem that holds not only for bounded but also for some unbounded
off-diagonal perturbations $V$.
\begin{introtheorem}
\label{Th1} Given a self-adjoint operator $A$ on the Hilbert space
$\fH$ assume that
\begin{equation*}
\spec(A) = \sigma_{-}\cup\sigma_{+} \text{\, and
\,}\sup\sigma_-<\inf\sigma_+.
\end{equation*}
Suppose that a symmetric operator $V$ on $\fH$ with
$\Dom(V)\supset\Dom(A)$ is off-diagonal with respect to the
decomposition $\fH=\Ran{\sE_A(\sigma_-)}\oplus\Ran{\sE_A(\sigma_+)}$
and the closure $L=\overline{A+V}$ of the sum $A+V$ with
$\Dom(A+V)=\Dom(A)$ is a self-adjoint operator. Then the spectrum of
$L$ consists of two subordinated components $\sigma'_-$ and
$\sigma'_+$ such that
\begin{equation*}
\sigma'_-\subset(-\infty,\sup\sigma_-], \quad
\sigma'_+\subset[\inf\sigma_+,+\infty),
\end{equation*}
and the following inequality holds
\begin{equation}
\label{Estin} \|\sE_{L}(\sigma'_-)-\sE_A(\sigma_-)\|\leq
\sin\left(\frac{1}{2}\arctan\varkappa\right),
\end{equation}
where
\begin{equation*}
%\label{jvIntro}
\varkappa= \inf_{\sup \sigma_{-} < \mu < \inf \sigma_{+}}
\sup\limits_{\mbox{\scriptsize$\begin{array}{cc} x\in \dom(A)\\
\|x\|=1\end{array}$}} \frac{|\lal x,J V x\ral|}{\lal x,|A-\mu|x\ral}
\end{equation*}
with $J$ given by \eqref{Jintro}.
\end{introtheorem}

Notice that throughout the paper we adopt the natural convention
that
$$
\arctan(+\infty)=\pi/2.
$$
In particular, under this convention
inequality \eqref{Estin} for $\varkappa=+\infty$ reads
$$\|\sE_{L}(\sigma'_-)-\sE_A(\sigma_-)\|\leq\frac{\sqrt{2}}{2}.$$

By Remark \ref{remTh1} $(iii)$ below the estimate \eqref{Estin} is
sharp.

Theorem \ref{Th1} is a corollary to a more general statement
(Theorem \ref{dkn}) that is valid even in the case where
$\sup\sigma_-=\inf\sigma_+$. In its turn, the Davis-Kahan $\tan
2\Theta$ Theorem (Theorem \ref{DK2T}) appears to be a simple
corollary to Theorem \ref{Th1}.

We also remark that for a class of unbounded off-diagonal
perturbations studied in \cite{AdLMSr} (cf. \cite{HMM},
\cite{MenShk}) the rough estimate
$\|\sE_{L}(\sigma'_-)-\sE_A(\sigma_-)\|\leq\frac{\sqrt{2}}{2}$ can
be proven by combining \cite[Theorem 5.3]{AdLMSr} and \cite[Theorem
5.6]{KMM5}. Example \ref{ksharp} to Theorem \ref{Th1} shows that
estimate \eqref{Estin} may hold (even with finite $\varkappa$) for
unbounded perturbations that do not fit the assumptions of
\cite{AdLMSr}.

As regards the spectral disposition \eqref{tant}, it has been proven
in \cite{KMM4} (see also \cite{KMM3}) that the gaps between
$\sigma_-$ and $\sigma_+$ remain open and the bound \eqref{edif}
holds if the perturbation  $V$ satisfies condition
\begin{equation*}
%\label{sqrt2}
\| V \| < \sqrt{2} d.
\end{equation*}
Moreover, under this condition by \cite[Theorems 1 (i) and
3.2]{KMM3} the following inclusions hold:
\begin{align}
\label{incl2} \sigma_{+}' &\subset
\mathbb{R}\backslash\Delta\text{\, and \,}
\sigma_{-}'\subset[\inf\sigma_--\delta_-,\sup\sigma_-+\delta_+],
\end{align}
where $\Delta=(\alpha,\beta)$, $\alpha<\beta$, stands for the finite
gap in the set $\sigma_+$ that contains $\sigma_-$ and
\begin{align}
\label{incl2a} \delta_-&=\|V\|\tan\left(\frac{1}{2}\arctan
\frac{2\|V\|}{\beta-\inf\sigma_-}\right)<\inf\sigma_- -\alpha,\\
\label{incl2b} \delta_+&=\|V\|\tan\left(\frac{1}{2}\arctan
\frac{2\|V\|}{\sup\sigma_- -\alpha}\right)<\beta-\sup\sigma_-.
\end{align}

The only known sharp bound  \cite[Theorem 2.4]{KMM4} (see also
\cite[Theorem 2]{KMM3}) for the norm of the difference
$\sE_{A+V}(\sigma'_{-})-\sE_{A}(\sigma_{-})$ involves the distance
from the initial spectral set $\sigma_+$ to the perturbed spectral
set $\sigma'_{-}$, and thus this bound is an \textit{a posteriori}
estimate.

Our second principal result just adds an \textit{a priori} sharp
bound for the norm
{$\|\sE_{A+V}(\sigma'_{-})-\sE_{A}(\sigma_{-})\|$} in the case where
\eqref{tant} holds and $\|V\|<d$.

\begin{introtheorem}
\label{main0} Given a self-adjoint operator $A$ on the Hilbert space
$\got H$ assume that
\begin{equation*}
\spec(A) = \sigma_{-}\cup\sigma_{+},  \quad
\dist(\sigma_{+},\sigma_{-}) = d > 0, \quad  \text{and} \quad
\sigma_{+} \cap\mathop{\mathrm{conv}}( \sigma_{-} )   = \emptyset.
\end{equation*}
Let $V$ be a bounded self-adjoint operator on $\fH$ off-diagonal
with respect to the decomposition $\got H = \ran
\sE_{A}(\sigma_{-})\oplus\ran \sE_A(\sigma_{+})$. Assume in addition
that
\begin{equation}
\label{Vlessd} \| V \| < d.
\end{equation}
Then
\begin{equation}
\label{edif1} \| \sE_{L}(\sigma'_{-}) - \sE_{A}(\sigma_{-})\| \le
\sin \Bigl( \arctan \frac{\|V\|}{d}\Bigr) = \frac{\|V\|}{\sqrt{d^2 +
\|V\|^2}},
\end{equation}
where $L=A+V$ with $\Dom(L)=\Dom(A)$.
\end{introtheorem}

\begin{introremark}
\label{raTanT} Estimate \eqref{edif1} can be equivalently written in
the form
\begin{equation}
\label{aTanT} \|\tan\Theta\|\leq\frac{\|V\|}{d},
\end{equation}
where $\Theta$ is the operator angle between the subspaces
$\Ran\sE_A(\sigma_-)$ and $\Ran\sE_{L}(\sigma'_-)$. Thus, Theorem
\ref{main0} may be called the \emph{a priori} $\tan\Theta$ Theorem.
It adds a new item to the list of fundamental estimates on the norm
of the difference of spectral projections known as $\sin \Theta$, $
\sin 2\Theta$, $\tan2\Theta$ Theorems (from \cite{Davis:123,DK}) and
\emph{a posteriori} $\tan\Theta$ Theorem (from \cite{DK,KMM3}).

\end{introremark}

We perform the proofs of both Theorems \ref{Th1} and \ref{main0} by
constructing the direct rotation \cite{Davis} from the subspace
$\ran \sE_{A}(\sigma_{-})$ to the subspace $\ran
\sE_{L}(\sigma_{-}')$.

Recall that the direct rotation $U$ from a closed subspace $\fM$ of
a Hilbert space $\fH$ to a closed subspace $\fN\subset\fH$ with
$\dim(\fM\cap \fN^{\perp})=\dim(\fM^\perp\cap \fN)$ is a unitary
operator on $\fH$ mapping $\fM$ onto $\fN$ and being such that for
any other unitary $W$ on $\fH$ with $\Ran W|_\fM=\fN$ the following
inequality holds: $\|I-U\|\leq\|I-W\|$ where $I$ is the identity
operator on $\fH$. That is, the direct rotation is closer (in the
operator norm topology) to the identity operator than any other
unitary operator on $\fH$ mapping $\fM$ onto $\fN$. The norm of the
difference between the corresponding orthogonal projections onto
$\fM$ and $\fN$ is completely determined by location of $\spec(U)$
on the unit circumference.

We extract information on the spectrum of the direct rotation from
the subspace $\ran\sE_{A}(\sigma_{-})$ to the subspace
$\ran\sE_{L}(\sigma_{-}')$ by using the following auxiliary result
which, we think, is of independent interest.

\begin{introtheorem}
\label{relemma} Let $T$ be a closed densely defined operator on a
Hilbert space $\got H$ with the polar decomposition $T = W|T|$.
Assume that $G$ is a bounded operator on $\got H$ such that both
$GT$ and $G^*T^*$ are accretive (resp. strictly accretive). Then the
products $GW$ and $WG$ are also accretive (resp. strictly accretive)
operators.
\end{introtheorem}
Notice that in this theorem and below an operator $T$ on the Hilbert
space $\got H$ is called \emph{accretive} (resp. \emph{strictly
accretive}) if
\begin{equation*}
%\label{accr}
\re\lal x,Tx\ral\geq 0 \text{\, (resp. \,}\re\lal x,Tx\ral> 0)
\text{\, for any\,} x\in\dom(T), \|x\|=1.
\end{equation*}
We also adopt the convention that the partial isometry $W$ in the
polar decomposition $T = W|T|$ is extended to $\Ker(T)$ by
\begin{equation}
\label{conviso} W|_{\Ker(T)}=0.
\end{equation}
In this way the isometry $W$ is uniquely defined on the whole space
$\fH$ (see, e.g., \cite[\S VI.7.2]{Kato}).

A convenient way to construct the direct rotation between two closed
subspaces of a Hilbert space is rendered by using a pair of
self-adjoint involutions associated with these subspaces. Although
the relative geometry of two subspaces is studied in great detail
(see, e. g., \cite{Halmos:69}, \cite{Kato}, \cite{SzNagy}), for
convenience of the reader we give in Section \ref{Beginning} a short
but self-contained exposition of the subject reformulating some
results in terms of a pair of involutions.

The remaining part of the article is organized as follows. Section
\ref{SPolar} contains a proof of Theorem~\ref{relemma}. The
principal result of this section is Theorem \ref{inv} that allows
one to compare two involutions one of which is associated with a
self-adjoint operator. Theorem \ref{Th1} and several other related
statements are proven in Section \ref{Stan2t}. Section \ref{S3isl}
contains a proof of Theorem~\ref{main0}. Notice that
Theorem~\ref{main0} appears to be a corollary to a more general
statement (Theorem \ref{trio1}) proven under a weaker than
\eqref{Vlessd} but more detail assumption \eqref{Vdd} involving the
length of the finite gap in $\sigma_+$ that contains the other
spectral set $\sigma_-$.

We conclude the introduction with description of some more notations
used throughout the paper. The identity operator on any Hilbert
space $\fH$ is denoted by $I$. Given a linear operator $T$ on $\fH$,
by $\mathcal{W}(T)$ we denote its numerical range,
\begin{equation*}
\mathcal{W}(T)= \{\lambda\in\mathbb{C}\,|\,\lambda=\lal
x,Tx\ral\text{ for some }x\in\Dom(T),\|x\|=1\}.
\end{equation*}
We use the standard concepts of commuting and anticommuting
operators dealing only with the case where at least one of the
operators involved is bounded~(see, e.g., \cite[\S
3.1.1]{Birman-Solomjak}). Assuming that $S$ and $T$ are operators on
$\fH$ suppose that the operator $S$ is bounded. We say that the
operators $S$ and $T$ \textit{commute}  (resp. \textit{anticommute})
and write $S\smile{T}$ or $T\smile{S}$ (resp. $S\frown{T}$ or
$T\frown{S}$) if $ST\subset TS$ (resp. $ST\subset -TS$).

\section{A pair of involutions}
\label{Beginning}

\subsection{An involution}

We start with recalling the concept of a (self-adjoint) involution
on a Hilbert space. This concept is a main tool we use in the
present paper. Notice that in the theory of spaces with indefinite
metric the involutions are often called canonical symmetries (see,
e.g., \cite{AI}).

\begin{defn}
\label{definv} A linear operator $J$ on the Hilbert space $\got H$
is called an \textit{involution} if
\begin{equation}
\label{JJd} J^* = J \quad\text{and}\quad J^2 = I.
\end{equation}
\end{defn}

In particular, if $P^-$ and $P^+=I-P^-$ are two complementary
orthogonal projections on $\got H$  then the differences $P^{+} -
P^{-}$ and $P^{-}-P^{+}$ are involutions.

By definition, any involution $J$ is a self-adjoint operator. In
fact, it is also a unitary operator since \eqref{JJd} yields
$J^*=J^{-1}$. Hence $\spec(J) = \{-1,1\}$ and the spectral
decomposition of $J$ reads
\begin{equation*}
%\label{JEE}
J = \int_{\mathbb{R}} \lambda {\sf E}_{J}(d \lambda) = {\sf
E}_J(\{+1\}) - {\sf E}_J(\{-1\}),
\end{equation*}
which implies that any involution on $\fH$ is the difference between
two complementary orthogonal projections. Obviously, the projections
${\sf E}_J(\{\pm 1\})$ are equal to
\begin{equation}
\label{Pinv} {\sf E}_J(\{+1\}) = \frac{1}{2} (I + J)
\quad\text{and}\quad {\sf E}_J(\{-1\}) = \frac{1}{2}  (I - J).
\end{equation}

\begin{defn}
Let $J$ be an involution on the Hilbert space $\got H$. The
subspaces
\begin{equation}
\label{Hpm} \fH_-= \Ran{\sf E}_J(\{-1\}) \quad\text{and}\quad \fH_+=
\Ran{\sf E}_J(\{+1\})
\end{equation}
are called the \textit{negative} and \textit{positive} subspaces of
the involution $J$, respectively. The decomposition
\begin{equation}
\label{Hsum} \fH=\fH_-\oplus\fH_+
\end{equation}
of $\got H$ into the orthogonal sum of the subspaces \eqref{Hpm} is
said to be \emph{associated} with $J$.
\end{defn}

Recall that a linear operator $A$ on $\fH$ is called diagonal with
respect to  decomposition~\eqref{Hsum} if the subspace $\fH_-$ (and
hence the subspace $\fH_+$) reduces $A$. A linear operator $V$ on
$\fH$ is said to be off-diagonal with respect to
decomposition~\eqref{Hsum} if
\begin{equation*}
\fH_-\cap\dom(V) = \ran P^{-}|_{\dom(V)}, \quad \fH_+\cap\dom(V) =
\ran P^{+}|_{\dom(V)},
\end{equation*}
where $P^-$ and $P^+$ are orthogonal projections onto $\got H_{-}$
and $\got H_{+}$, respectively, and
\begin{equation}
\label{Voff} \ran V|_{\fH_-\cap\dom(V)} \subset\fH_{+}, \quad \ran
V|_{\fH_+\cap\dom(V)} \subset\fH_{-}.
\end{equation}

A criterion for an operator on $\got H$ to be diagonal or
off-diagonal with respect to the orthogonal decomposition of $\got
H$ associated with an involution $J$ can be formulated in terms of a
commutation relation between this operator and $J$.

\begin{lem}
A linear operator $A$ on the Hilbert space $\got H$ is diagonal with
respect to the orthogonal decomposition of $\got H$ associated with
an involution $J$ if and only if $J\smile A$.
\end{lem}
\begin{proof}
This assertion is an immediate corollary to \cite[Theorem 1 in \S
3.6]{Birman-Solomjak}.
\end{proof}

\begin{lem}
A linear operator $V$ on the Hilbert space $\got H$ is off-diagonal
with respect to the orthogonal decomposition of $\got H$ associated
with an involution $J$ if and only if $J\frown V$.
\end{lem}

\begin{proof} \textit{``Only if part.''}
Assume that $V$ is off-diagonal with respect to an orthogonal
decomposition of $\got H$ associated with $J$. Let $P^{\pm} =
\sE_{J}(\{\pm 1\})$. Then $J = P^{+} - P^{-}$ and $P^{+} + P^{-} =
I$. By the hypothesis one infers that $P^\pm x\in\Dom(V)$  for any
$x\in\Dom(V)$. Hence $x\in\dom(V)$ implies $Jx\in\Dom(V)$. Moreover,
for any $x\in\dom(V)$ the following chain of equalities holds
\begin{align*}
VJx & =VP^+x-VP^-x\\
    & =P^-VP^+x-P^+VP^-x\\
    & =P^-V(P^++P^-)x-P^+V(P^++P^-)x\\
    & =(P^--P^+)Vx \\
    & = - JV x,
\end{align*}
since $P^+ V P^+x = P^-VP^-x=0$ (cf.~\eqref{Voff}). Thus $J \frown
V$.

\textit{``If part.''} Suppose that $J\frown V$ which means that
$(i)$ $x\in \Dom(V)$ implies $J x\in \dom(V)$ and $(ii)$ $VJx =
-JVx$ for all $x\in \dom(V)$. Let $\got H_{\pm} = \ran
\sE_{J}(\{\pm\})$. Condition $(i)$ and equalities \eqref{Pinv} imply
that $\sE_{J}(\{\pm 1\}) x \in \dom(V)$ whenever $x\in \dom(V)$.
Therefore it follows from condition $(ii)$ that if $x_{-} \in\got
H_{-} \cap\dom(V) $, then $V x_{-} = -V J x_{-} = J V x_{-}$. Hence
$V x_{-} \in \got H_{+}$ for all $x_{-} \in \got H_{-}\cap \dom(V)$.
In a similar way one verifies that $V x_{+} \in \got H_{-}$ for all
$x_{+} \in \got H_{+} \cap\dom(V)$. Hence $V$ is off-diagonal with
respect to the decomposition of $\got H$ associated with $J$, which
completes the proof.
\end{proof}

\begin{rem}
\label{rematr} Operators that are diagonal or off-diagonal with
respect to the decomposition~\eqref{Hsum} are often written in the
block operator matrix form,
\begin{equation*}
%\label{AVmatr}
A=\left(\begin{array}{ll} A_- & 0 \\ 0 & A_+
\end{array}\right), \quad
V=\left(\begin{array}{ll} 0 & V_+ \\ V_- & 0
\end{array}\right),
\end{equation*}
where $A_\pm$ are the parts of the diagonal operator $A$ in $\got
H_{\pm}$, and $V_\pm$ are the corresponding restrictions of the
off-diagonal operator $V$ to $\got H_{\pm}$,
\begin{equation*}
%\label{Vparts}
A_{\pm} = A|_{\dom(A)\cap\fH_\pm}, \quad
V_\pm=V|_{\dom(V)\cap\fH_\pm}.
\end{equation*}

In particular, if both $A$ and $V$ are closed operators and, in
addition, $V$ is bounded, then the closed operator $L=A+V$ with
$\dom(L)=\dom(A)$ admits the block operator matrix representation
\begin{equation}
\label{Lmatr} L=\left(\begin{array}{ll} A_- & V_{+} \\ V_{-} & A_+
\end{array}\right).
\end{equation}
In this case
\begin{equation*}
A = \frac{1}{2} ( L + J L J), \quad V = \frac{1}{2} \overline{( L -
J L J )},
\end{equation*}
where $J$ is the involution that corresponds to the
decomposition~\eqref{Hsum}.
\end{rem}

Notice that the study of invariant subspaces for block operator
matrices of the form \eqref{Lmatr} is closely related to the
question concerning existence of solutions to the associated
operator Riccati equations (see, e.g., \cite{AMM} and references
therein).

\subsection{Involutions in the acute case}

Recall that two closed subspaces $\fM$ and $\fN$ of a Hilbert space
$\got H$ are said to be in the acute case if
\begin{equation*}
\fM\cap\fN^\perp = \{0\} \quad\text{and}\quad \fM^\perp\cap\fN =
\{0\}.
\end{equation*}
To formulate the notion of the acute case in terms of the
corresponding involutions we adopt the following definition.

\begin{defn}
\label{defacute} Involutions $J$ and $J'$ on the Hilbert space $\got
H$ are said to be in the \textit{acute case} if
\begin{equation*}
\ker(I + J'J) = \{0\}.
\end{equation*}
\end{defn}

\begin{rem}
\label{remker} By inspection, $\ker(I + J' J) = \ker(I + J J')$
which means that this definition is symmetric with respect to the
entries $J$ and $J'$.
\end{rem}

\begin{lem}
\label{remcom} If involutions $J$ and $J'$ are in the acute case and
$J\smile J'$, then $J=J'$.
\end{lem}
\begin{proof}
Taking into account the self-adjointness of both $J$ and $J'$, the
hypothesis $JJ'=J'J$ implies that the unitary operator $J'J$ is
self-adjoint. Hence $\spec(J'J)\subset\{-1,1\}$. Then from the
assumption that $J$ and $J'$ are in the acute case it follows that
$-1\not\in\spec(J'J)$. This yields $J'J = I$ and hence $J=J'$.
\end{proof}

Some criteria for a pair of involutions $J$ and $J'$ to be in the
acute case are presented in Lemma \ref{acute} below. In particular,
this lemma justifies Definition~\ref{defacute} stating that $J$ and
$J'$ are in the acute case if and only if their negative (resp.
positive) subspaces are in the acute case.

One of the criteria in Lemma \ref{acute} involves the numerical
range $\cW(J'J)$ of the product $J'J$. Since $J'J$ is a unitary
operator, its numerical range is a subset of the unit disc
$\{\lambda\in\mathbb{C}\,|\,\,|\lambda|\le 1\}$. Equalities $J' J =
J(JJ')J=J(JJ')J^{-1}$ imply that the products $J'J$ and $JJ'$ are
unitarily equivalent. Hence $\cW(J'J) = \cW(JJ')$. By $J J' =
(J'J)^*$ this means that the numerical range of $J'J$ is symmetric
with respect to the real axis.

\begin{lem}
\label{acute} Let $J$ and $J'$ be two involutions on the Hilbert
space $\got H$. Assume that $\got H_{\pm} = \ran \sE_{J}(\{\pm 1\})$
and $\got H_{\pm}' = \ran\sE_{J'}(\{\pm 1\})$. The following four
statements are equivalent:
\begin{itemize}

\item[$(i)$] $\quad\got H_-\cap\got H'_+ = \{0\}
\quad\text{and}\quad \got H_+\cap\got H'_- = \{0\}$,

\smallskip

\item[$(ii)$] $\quad\ker( I + J' J ) = \{0\}$,

\smallskip

\item[$(iii)$] $\quad \|(J'-J) x \| < 2 \| x \|
\quad \text{for all}\quad x \in \got H,\quad x\ne 0$,

\smallskip

\item[$(iv)$] $\quad-1\not\in \cW(J'J)$.
\end{itemize}
\end{lem}

\begin{proof}
We prove the implications $(i)\Rightarrow (ii) \Rightarrow (iii)
\Rightarrow (iv) \Rightarrow (i)$.

$(i) \Rightarrow (ii)$. We prove this implication by contradiction.
Suppose that $\ker (I + J' J)\neq\{0\}$ and $x\in\ker(I+ J'J)$ is a
non-zero vector. Representing this vector as $x = x_- + x_+$ with
$x_-\in \got H_-$ and $x_+ \in \got H_+$ one obtains $(I + J' J) x =
(I - J') x_- + (I + J') x_+$ and hence
\begin{equation}
\label{sum} (I - J') x_- + (I + J') x_+=0
\end{equation}
since $(I + J' J)x=0$. Applying $(I - J')$ to both parts of
\eqref{sum} gives $(I - J')^2 x_- = 0$ and thus $J' x_- = x_-$.
Therefore $x_-$ is an eigenvector of the operator $J'$ corresponding
to the eigenvalue $+1$ which means $x_-\in \got H_- \cap \got H'_+$.
In a similar way, by applying  $(I + J')$ to both parts of
(\ref{sum}), one concludes that $J' x_+ = - x_+$ and hence $x_+ \in
\got H_+ \cap \got H'_-$. Then it follows from condition $(i)$ that
$x_- = x_+ = 0$ and thus $x=0$ which contradicts the assumption.

$(ii) \Rightarrow (iii)$. It follows from condition $(ii)$ that
$\|(I + J' J)x\|>0$ for any non-zero $x\in \got H$. Then by taking
into account the identities
\begin{equation*}
\|(J - J') x\|^2+\|(J + J') x \|^2  = 4 \| x \|^2
\end{equation*}
and
\begin{equation*}
\|(J+J')x\| = \| J'(J'+J)x\| = \|(I + J' J)x\|
\end{equation*}
one easily concludes that $(ii)$ implies $(iii)$.

$(iii) \Rightarrow (iv)$. By inspection
\begin{equation*}
\| x \|^2 + \re \lal x, J'J x\ral = \frac{1}{2} \Bigl\{ 4 \| x\|^2 -
\|(J-J')x\|^2 \Bigr\}.
\end{equation*}
Hence $(iii)$ implies
\begin{equation*}
\| x \|^2 + \re \lal x, J'J x\ral>0 \quad\text{for any non-zero
$x\in \got H$}.
\end{equation*}
In particular, this means that $\re\lal x, J'J x\ral>-1$ for any
$x\in\fH$ such that $\|x\|=1$ and therefore $-1\not\in \cW(J'J)$.

$(iv)\Rightarrow (i)$. Suppose that at least one of the subspaces
$\got H_-\cap\got H'_+$ and $\got H_+ \cap \got H'_-$ is
non-trivial. Pick up vectors $x_-\in \got H_-\cap\got H'_+$ and $x_+
\in \got H_+ \cap \got H'_-$ such that at least one of them is
non-zero. Clearly, $J' J (x_- + x_+) = J' (-x_- + x_+) = -(x_- +
x_+)$ which means that $-1$ is an eigenvalue of the operator $J' J$
and thus $-1\in \cW(J'J)$. This contradicts the assumption $(iv)$
and thus proves the implication.
\end{proof}

\begin{rem}
\label{remdif} Making use of relationship~\eqref{Pinv} between an
involution and its spectral projections yields
\begin{equation*}
P'^{+} - P^{+} = P^{-} - P'^{-} =\frac{J'-J}{2},
\end{equation*}
where $P^\pm={\sf E}_{J}(\{\pm1\})$ and $P'^\pm={\sf
E}_{J'}(\{\pm1\})$.
\end{rem}

\begin{cor}
If
\begin{equation*}
\|P'^{-} - P^{-}\| < 1 \quad (\text{or  \, \,} \|P'^{+} - P^{+}\| <
1)
\end{equation*}
holds then the involutions $J$ and $J'$ are in the acute case.
Hence, the negative (resp. positive) subspaces of $J$ and $J'$ are
also in the acute case.
\end{cor}

\subsection{The direct rotation}
\label{secDiRot} Let $J$ and $J'$ be involutions on $\got H$. Assume
that $\got H_{-}$ and $\got H_{+}$ are the negative and positive
subspaces of $J$, respectively. Similarly, assume that $\got H_{-}'$
and $\got H_{+}'$ are the negative and positive subspaces of $J'$.
It is well known (see, e.g., \cite[Theorem 3.1]{Davis}) that if
\begin{equation}
\label{dims} \dim( \got H_{-}\cap\got H'_{+}) = \dim (\got
H_{+}\cap\got H'_{-}),
\end{equation}
then there exists a unitary operator $W$ on $\got H$ mapping $\got
H_{-}$ onto $\got H_{-}'$ and $\got H_{+}$ onto $\got H_{+}'$.
Clearly, $W$ satisfies the commutation relation
\begin{equation}
\label{wcomm} J' W = W J.
\end{equation}
In particular, by Lemma~\ref{acute} such a unitary $W$ exists if $J$
and $J'$ are in the acute case. The canonical choice of the unitary
mapping of one subspace in the Hilbert space onto another, the
so-called direct rotation, was suggested by C. Davis in \cite{Davis}
and T. Kato in \cite[Sections I.4.6 and I.6.8]{Kato}. The idea of
this choice goes back yet to  B. Sz.-Nagy (see \cite[\S
105]{SzNagy}). We adopt the following definition of the direct
rotation.

\begin{defn}
\label{DiRot} Let $J$ and $J'$ be involutions on the Hilbert space
$\fH$. A unitary operator $U$ on $\fH$ is called {\it the direct
rotation} from $J$ to $J'$ if
\begin{equation}
\label{DiRotf} (i) \quad J' U = U J, \quad (ii) \quad U^2 = J'J,
\quad (iii) \quad \re U \ge 0.
\end{equation}
\end{defn}
\begin{rem}
The spectrum of any direct rotation is a subset of the unit
circumference lying in the closed right half-plane symmetrically
with respect to the real axis. To see this, observe that equalities
$(i)$ and $(ii)$ imply $U^* = JUJ$ by taking into account that $U$
is a unitary operator. Hence the operator $U$ is unitary equivalent
to its adjoint and thus the spectrum of $U$ is symmetric with
respect to the real axis. From $(iii)$ it follows that this spectrum
is a subset of the half-plane $\{z\in\bbC\,|\,\re z\geq0\}$. To
complete the proof of the statement it only remains to recall that
the spectrum of any unitary operator lies on the unit circumference.
\end{rem}

We give a short proof of the existence and uniqueness of the direct
rotation for the instance where the corresponding involutions are in
the acute case. For a different proof of this fact see
\cite[Propositions 3.1 and 3.3]{DK}.

\begin{thm}
\label{rot} If involutions $J$ and $ J'$ are in the acute case then
there is a unique direct rotation from $J$ to $J'$.
\end{thm}
\begin{proof}
We divide the proof into two parts. In the first part we prove the
existence of a direct rotation from $J$ to $J'$. The uniqueness of
the direct rotation is proven in the second part.
\medskip

\noindent \textit{(Existence.)} Set $ T = I + J' J. $ One easily
verifies that $T$ is a normal operator. By hypothesis
\begin{equation}
\label{kerTZ} \ker(T) = \ker(T^*) = \{0\}
\end{equation}
taking into account Remark \ref{remker}. Hence the the isometry $U$
in the polar decomposition
\begin{equation}
\label{polarTU} T = U |T| = |T| U,
\end{equation}
is a unitary operator (see \cite[\S 110]{SzNagy}).

By inspection
\begin{equation}
\label{commT} J' T = T J
\end{equation}
and thus
\begin{gather*}
J |T|^2 = J T^* T = T^* J' T = T^* T J = |T|^2 J,\\
J' |T|^2 = J' T T^* = T J T^* = T T^* J' = |T|^2 J'.
\end{gather*}
Hence $J \smile |T|$ and $J'\smile |T|$. Then \eqref{polarTU} and
\eqref{commT} yield $|T|(J' U - U J) = 0$, which implies that
\begin{equation}
\label{JUUJ} J' U = U J
\end{equation}
since  $\ker(|T|)=\ker (T) = \{0\}$. Observing that $J'J T^* = T$,
by the same reasoning one obtains $|T|( U - J' J U^* ) =0$. Hence $U
= J' J U^*$ and thus
\begin{equation}
\label{U2JJ} U^2 = J' J.
\end{equation}
Finally, $T + T^* = |T|^2$ and $T + T^* = |T|( U + U^*)$ imply
$|T|(U+U^* - |T|) = 0$. Therefore
\begin{equation}
\label{ReU} \re U =\frac{1}{2} |T| \ge 0.
\end{equation}
Comparing \eqref{JUUJ}, \eqref{U2JJ}, and \eqref{ReU} with
\eqref{DiRotf}, one concludes that $U$ is the direct rotation from
$J$ to $J'$.
\smallskip

\noindent\textit{(Uniqueness.)} Suppose that $ U'$ is another
unitary operator such that $ U'^2 = U^2$ and $\re  U' \ge 0$. By
inspection,
\begin{align*}
(\re  U')^2 & = \frac{1}{2} \Bigl(I+ \re(U'^2) \Bigr)
 =\frac{1}{2} \Bigl(I+\re(U^2)\Bigr)
 = (\re  U)^2.
\end{align*}
Then it follows from the uniqueness of the positive square root of a
positive operator that $\re  U = \re U'$. In addition, the
requirement $\im ( U^2 ) = \im (U'^2)$ implies $\re  U (\im  U - \im
U') =0$ which means that $\im U = \im U'$ since $\ker(\re U) =
\ker(|T|)=\{0\}$ by combining \eqref{kerTZ} and \eqref{ReU}. Thus
$U' = \re U + i \im U = U$, completing the proof.
\end{proof}

\begin{rem}
In the nonacute case the direct rotation exists if and only
if~\eqref{dims} holds (see \cite[Proposition 3.2]{DK}). If it
exists, it is not unique.
\end{rem}

To specify location of the spectrum of a unitary operator on the
unit circumference we introduce the notion of the spectral angle.

\begin{defn}
\label{spang} Let $W$ be a unitary operator. The number
\begin{equation*}
\vartheta(W) = \sup_{z\in \spec(W)} |\arg z|, \quad \arg
z\in(-\pi,\pi],
\end{equation*}
is called the \textit{spectral angle} of $W$.
\end{defn}

\begin{rem}
\label{wwad} $\vartheta(W^*)=\vartheta(W).$
\end{rem}

\begin{rem}
The (self-adjoint) operator angle between two closed subspaces in a
Hilbert space is expressed through the direct rotation $U$ from one
of these subspaces to the other one by $\Theta=\arccos(\re U)$ (see
\cite[Eq. (1.18)]{DK}). Hence $\vartheta(U)$ is nothing but the
spectral radius of the corresponding operator angle $\Theta$.
\end{rem}

The next statement shows that the spectral angle $\vartheta(W)$ is a
quantity that characterizes the distinction of the unitary operator
$W$ from the identity operator.

\begin{lem}
\label{spanlem} Let $W$ be a unitary operator. Then
\begin{equation}
\label{ImW} \| I - W \| = 2 \sin\left(\frac{\vartheta(W)}{2}\right).
\end{equation}
\end{lem}

\begin{proof}
Observe that $I - W$ is a normal operator. Then by using the
spectral mapping theorem one concludes that the following chain of
equalities holds:
\begin{align*}
\| I - W \| & = \sup_{\lambda\in \spec(I - W)}|\lambda|\\
&= \sup_{z\in\spec(W)}|1-z|\\
&= \sup_{z\in\spec(W)} 2 \sin \left(\frac{|\arg z|}{2}\right)\\
&= 2 \sin\Bigl( \frac{1}{2} \sup_{z\in\spec(W)}|\arg z|\Bigr) \\
& = 2\sin \left(\frac{\vartheta(W)}{2}\right),
\end{align*}
where $\arg z \in (-\pi,\pi]$.
\end{proof}

\begin{rem}
If $U$ is the direct rotation from an involution $J$ to an
involution $J'$ then it possesses the extremal property
\begin{equation*}
\vartheta(U) \le \vartheta(W),
\end{equation*}
where $W$ is any other unitary operator satisfying \eqref{wcomm}.
This can be easily seen from \eqref{ImW} by using \cite[Theorem
7.1]{Davis} which states that $\|I-U\| \le \|I-W\|$.
\end{rem}

\begin{rem}
\label{star} Again assume that $U$ is the direct rotation from an
involution $J$ to an involution  $J'$. Then by \eqref{DiRotf} the
spectral mapping theorem implies
\begin{equation}
\label{uhalf} 0\leq\vartheta(U)\leq\frac{\pi}{2} \quad
\text{and}\quad \vartheta(U) = \frac{1}{2} \vartheta(J'J).
\end{equation}
Since $\| J' - J \| = \| I - J'J\|$, by \eqref{ImW} it follows from
\eqref{uhalf} that
\begin{equation*}
\| J' - J \| = 2 \sin \left(\frac{\vartheta(J'J)}{2}\right) = 2 \sin
\vartheta(U).
\end{equation*}
Hence by Remark~\ref{remdif}
\begin{equation}
\label{difPt} \| P'^{+} - P^{+} \|= \| P'^{-} - P^{-} \|=\sin
\vartheta(U),
\end{equation}
where $P^{\pm} = \sE_{J}(\{\pm 1\})$ and $P'^{\pm} = \sE_{J'}(\{\pm
1\})$.
\end{rem}

In the proof of the next lemma we will use the following notation.
Assume that $\mathcal{S}$ is a subset of the complex plane. Then
$e^{i\varphi}\mathcal{S}$ denotes the result of rotation of
$\mathcal{S}$ by the angle $\varphi\subset(-\pi,\pi]$ around the
origin, that is,
\begin{equation*}
e^{i\varphi}\mathcal{S}=\{z\in\bbC\,|\,\,z=e^{i \varphi}\zeta
\text{\, for some\,}\zeta\in\mathcal{S}\}.
\end{equation*}

\begin{lem}
\label{UU} Let  $W_1$ and $W_2$ be two unitary operators on the
Hilbert space $\got H$. Then
\begin{equation}
\label{sumU} |\vartheta(W_1) - \vartheta(W_2)| \le \vartheta(W_2
W_1) \le \vartheta(W_1) + \vartheta(W_2).
\end{equation}
\end{lem}

\begin{proof}
First, we prove inequality
\begin{equation}
\label{ttsum} \vartheta(W_2 W_1) \le \vartheta(W_1) + \vartheta(W_2)
\end{equation}

Denote by $\vartheta_1$, $\vartheta_2$ and $\vartheta_3$ the
spectral angles of $W_1$, $W_2$, and $W_2 W_1$, respectively. The
case $\vartheta_1 + \vartheta_2 \ge \pi$ is trivial since
$\vartheta_3 \leq \pi$ by Definition~\ref{spang}. If $\vartheta_1 +
\vartheta_2 < \pi$, we prove \eqref{ttsum} by contradiction. Suppose
that the opposite inequality holds, that is,
\begin{equation*}
%\label{trit}
\vartheta_3 > \vartheta_1 + \vartheta_2.
\end{equation*}
Then there is a number $\varphi\in (-\pi,\pi]$ such that
$e^{i\varphi}\in \spec(W_2 W_1)$ and
\begin{equation}
\label{phit} \vartheta_1+\vartheta_2 < |\varphi|\le \pi.
\end{equation}
Since $W_2 W_1$ is a normal (unitary) operator, there exists a
sequence of vectors $x_n\in\got H$, $n=1,2,...,$ such that
\begin{equation}
\label{xn} \|x_n\|=1 \text{\, and \,} \| W_2 W_1 x_n - e^{i\varphi}
x_n \| \to 0, \quad n\to \infty.
\end{equation}
Indeed, if $e^{i\varphi}$ is an eigenvalue of $W_2 W_1$, to satisfy
\eqref{xn} one simply takes $x_n=x_\varphi$, $n=1,2,\ldots$, where
$x_\varphi$ is a normalized eigenvector of $W_2 W_1$ corresponding
to the eigenvalue $e^{i\varphi}$, i.e. $W_2 W_1
x_\varphi=e^{i\varphi}x_{\varphi}$. Otherwise such a sequence exists
by the Weyl criterion for the essential spectrum.

Let $z_{1,n} = \lal x_n, W_1 x_n\ral$ and $z_{2,n} = \lal x_n, W_2^*
x_n \ral$. Clearly, \eqref{xn} yields
\begin{equation}
\label{lim1} | z_{1,n} - e^{i\varphi} z_{2,n} | \to 0, \quad n\to
\infty,
\end{equation}
since by the Schwartz inequality
\begin{align*}
| z_{1,n} - e^{i\varphi} z_{2,n} | = &
|\lal x_n, W_1 x_n - e^{i\varphi} W_2^* x_n\ral|\\
&\le \| W_1 x_n - e^{i\varphi} W_2^* x_n \| = \|W_2 W_1 x_n -
e^{i\varphi} x_n\|.
\end{align*}
Taking into account that $z_{1,n}\subset \cW(W_1)$ and
$z_{2,n}\subset \cW(W_2^*)$, from \eqref{lim1} one concludes that
\begin{equation}
\label{WW0} \dist\bigl(\cW(W_1),e^{i\varphi}\cW(W_2^*)\bigr)=0.
\end{equation}

Meanwhile, if $W$ is a unitary operator with the spectral angle
$\vartheta$, the spectral theorem implies
\begin{equation*}
%\label{incl}
\cW(W) \subset \mathcal{S}_\vartheta \quad\text{and}\quad \cW(W^*)
\subset \mathcal{S}_\vartheta,
\end{equation*}
where
\begin{equation*}
\mathcal{S}_\vartheta = \{ z\in \mathbb{C}\, | \,\, \re z \ge
\cos\vartheta \text{\, and \,} |z|\le 1   \}
\end{equation*}
is a segment of the closed unit disc centered at the origin.
Therefore, $\cW(W_1)\subset\mathcal{S}_{\vartheta_1}$ and
$\cW(W_2^*)\subset\mathcal{S}_{\vartheta_2}$. Obviously,
$e^{i\varphi}\cW(W_2^*)\subset
e^{i\varphi}\mathcal{S}_{\vartheta_2}$ and hence
\begin{equation}
\label{dWW} \dist\bigl(\cW(W_1),e^{i\varphi}\cW(W_2^*)\bigr)\geq
\dist\bigl(\mathcal{S}_{\vartheta_1},e^{i\varphi}
\mathcal{S}_{\vartheta_2}\bigr).
\end{equation}
One easily verifies by inspection that under the assumption
\eqref{phit}
\begin{equation*}
  \dist\bigl(\mathcal{S}_{\vartheta_1},e^{i\varphi}
\mathcal{S}_{\vartheta_2}\bigr) =
     2 \sin\left(\frac{|\varphi|-\vartheta_1-\vartheta_2}{2}\right)
    \sin\left(\frac{|\varphi|+\vartheta_2-\vartheta_1}{2}\right) > 0
\end{equation*}
and thus by \eqref{dWW}
\begin{equation*}
\dist\bigl(\cW(W_1),e^{i\varphi}\cW(W_2^*)\bigr)>0
\end{equation*}
which contradicts \eqref{WW0}. This completes the proof of
\eqref{ttsum}.

By Remark \ref{wwad}, inequality \eqref{ttsum} implies
\begin{gather}
\label{ineq1} \vartheta(W_2) = \vartheta(W_2 W_1 W_1^*) \le
\vartheta(W_2 W_1) +
\vartheta(W_1^*) = \vartheta(W_2 W_1) + \vartheta(W_1), \\
\label{ineq2} \vartheta(W_1) = \vartheta(W_2^* W_2 W_1) \le
\vartheta(W_2^*) + \vartheta(W_2 W_1) = \vartheta(W_2) +
\vartheta(W_2 W_1).
\end{gather}
Combining~\eqref{ineq1} and~\eqref{ineq2} yields the left inequality
in \eqref{sumU}. The proof is complete.
\end{proof}

\begin{rem}
Setting $W_1 = e^{i \vartheta_1} I$ and $W_2 = e^{i\vartheta_2}I$
with $\vartheta_1$, $\vartheta_2$ appropriate reals, one verifies
that both inequalities of~\eqref{sumU} are sharp.
\end{rem}

\section{A property of the polar decomposition}
\label{SPolar} In this section we give a proof of Theorem
\ref{relemma}. We also derive corollaries to this theorem for the
case where one of the operators involved is self-adjoint and the
other one is related to an involution.

We start with an auxiliary result.

\begin{lem}
\label{unit} Let $A$ be a positive operator on the Hilbert space
$\got H$. Suppose that $x,y\in\fH$ are such that
\begin{equation}
\label{rea1} \re \lal x, A(A^2 + \alpha)^{-1} y\ral > 0 \quad (\geq
0) \quad\text{for any}\quad\alpha>0.
\end{equation}
Then
\begin{equation}
\label{rea1xy} \re\lal x,Q y\ral > 0 \quad (\geq 0),
\end{equation}
where $Q$ is the orthogonal projection onto $\ker(A)^{\perp}$.
\end{lem}

\begin{proof}
By the spectral theorem
\begin{equation*}
\re \lal x, A(A^2 + \eta^2)^{-1} y\ral = \int_{\mathbb{R}}
\frac{\lambda m(d\lambda)}{\lambda^2 + \eta^2} = \int_{(0,+\infty)}
\frac{\lambda m(d\lambda)}{\lambda^2 + \eta^2}, \quad
0\neq\eta\in\bbR,
\end{equation*}
where for any Borel set $\delta\subset\bbR$ the Lebesgue--Stieltjes
measure $m(\delta)$ reads
\begin{equation*}
m(\delta)=\re\lal x,\sE_A\bigl(\delta)y\ral.
\end{equation*}
Hence for any $\varepsilon > 0$
\begin{align*}
\int\limits_{\varepsilon}^{1/\varepsilon} d\eta\,\re\lal x, A(A^2 +
\eta^2)^{-1} y\ral &= \int\limits_{\varepsilon}^{1/\varepsilon} d
\eta \int_{(0,+\infty)}
\frac{\lambda m(d\lambda)  }{\lambda^2 + \eta^2}   \\
& =  \int_{(0,+\infty)} m(d\lambda)\,
\int_{\varepsilon}^{1/\varepsilon} \frac{ \lambda d \eta
}{\lambda^2 + \eta^2}
\end{align*}
by the Fubini theorem. Therefore
\begin{align}
\label{deta1} \int\limits_{\varepsilon}^{1/\varepsilon}
d\eta\,\re\lal x, A(A^2 + \eta^2)^{-1} y\ral &=
\int_{(0,+\infty)}m(d\lambda)\,\,\Bigl[
\arctan\left(\frac{1}{\lambda\varepsilon}\right) -
\arctan\left(\frac{\varepsilon}{\lambda}\right)\Bigl].
\end{align}
From \eqref{deta1} one immediately infers that
\begin{equation}
\label{pim2} \lim\limits_{\varepsilon\downarrow 0}\,\,
\int\limits_{\varepsilon}^{1/\varepsilon} d\eta\,\re\lal x, A(A^2 +
\eta^2)^{-1} y\ral=\frac{\pi}{2}\, m\bigl((0,+\infty)\bigr).
\end{equation}
Notice that $m\bigl((0,+\infty)\bigr)=\re\lal x,Q y\ral$ since $Q =
\sE_{A}\bigl((0,+\infty)\bigr)$. Hence \eqref{pim2} yields
\begin{equation}
\label{lim} \re\lal x,Q y\ral=\lim_{\varepsilon\downarrow 0}\,\,
\frac{2}{\pi}\,\,\int\limits_{\varepsilon}^{1/\varepsilon} d \eta
\,\re \lal x,A(A^2 + \eta^2)^{-1} y\ral.
\end{equation}
Clearly, by \eqref{lim} inequalities \eqref{rea1xy} follow directly
from the corresponding assumptions \eqref{rea1}. The proof is
complete.
\end{proof}

With Lemma \ref{unit} we are ready to prove Theorem \ref{relemma}.

\begin{proof}[Proof of Theorem~\ref{relemma}]

Assume first that the operators $GT$ and $G^*T^*$ are both
accretive. To prove that $GW$ is also an accretive operator, pick up
arbitrary $\alpha>0$ and $x\in \got H$ and set
\begin{equation}
\label{gh} g = (T^* T + \alpha)^{-1} x.
\end{equation}
Taking into account that
%\begin{equation}
%\label{gh2}
$ g \in \dom(T), $
%\end{equation}
introduce
\begin{equation}
\label{gTh} h =T g= T(T^* T + \alpha)^{-1} x.
\end{equation}
Clearly, $h\in\dom(T^*)$ and
\begin{equation}
\label{x} x = \alpha g + T^* h.
\end{equation}
By using \eqref{gh}, \eqref{gTh}, and \eqref{x} it is easy to verify
that the following chain of equalities holds
\begin{align}
\nonumber \re\lal  W^* G^* x, |T| (|T|^2 + \alpha)^{-1} x\ral &=
\re\lal  G^* x, W |T| g \ral\\
\nonumber
&=\re\lal  G^* x, T g\ral \\
\nonumber
&=\re \lal x, G h\ral \\
\nonumber
&=\re\lal \alpha g + T^*h, G h\ral\\
\nonumber &=\alpha \re\lal g, G \nonumber
h\ral  + \re\lal  G h, T^* h\ral\\
\label{long} &=\alpha \re\lal  g, G T g \ral  + \re\lal h, G^* T^*
h\ral.
\end{align}
Since by hypothesis both $GT$ and $G^*T^*$  are accretive,
\eqref{long} implies that
\begin{equation*}
\re\lal  W^* G^* x, |T| (|T|^2 + \alpha)^{-1} x\ral \geq 0 \text{\,
for any \, $\alpha>0$ \, and \,} x\in\fH,
\end{equation*}
and hence by Lemma \ref{unit}
\begin{equation*}
\re\lal W^* G^* x,Q x\ral=\re\lal x, G W Q x\ral\geq0,
\end{equation*}
where $Q$ is the orthogonal projection onto $\ker(|T|)^{\perp}$.
According to the convention \eqref{conviso} we have $\ker(|T|) =
\ker(T) = \ker(W)$. Then one concludes that $W Q = W$ and hence
\begin{equation*}
\re\lal x, G W x \ral\geq0 \quad \text{for all}\quad x\in \got H,
\end{equation*}
which proves that the operator $GW$ is accretive.

Further, assume that $GT$ and $G^*T^*$ are both strictly accretive
operators. In particular, this implies that
\begin{equation}
\label{kerT0} \ker(T) = \ker(|T|) =\{0\}.
\end{equation}
In this case if $x\neq0$ then neither $g$ nor $h$ defined in
\eqref{gh} and \eqref{gTh} can be zero vectors. Indeed, the equality
$g= 0$ implies $h=Tg=0$ and hence by \eqref{x} it contradicts the
assumption $x\ne 0$. Independently, the equality $h=0$ yields $g\in
\ker(T)$ by taking into account \eqref{gTh}. Then $x\in \ker(T)$
since $x = \alpha g$ by \eqref{x}. This is again a contradiction
because of \eqref{kerT0}.

Therefore if $x\neq 0$ and $\alpha>0$ then necessarily $g\neq0$,
$h\neq0$. Hence by \eqref{long} now we have the strict inequality
\begin{equation*}
%\label{GUG}
\re\lal  W^* G^* x, |T| (|T|^2 + \alpha)^{-1} x\ral > 0.
\end{equation*}
Then by taking into account \eqref{kerT0} Lemma~\ref{unit} proves
the strict accretiveness of the operator~$GW$.

The accretiveness (resp., the strict accretiveness) of the operator
$WG$ can be proven in a similar way.
\end{proof}

Now assume that $T$ is a self-adjoint operator on the Hilbert space
$\fH$ and $\ker(T)=\{0\}$. Then the isometry $J'$ in the polar
decomposition
\begin{equation}
\label{polarT} T = J' |T|
\end{equation}
is an involution that reads
\begin{equation*}
%\label{difJL}
J' = \sE_{T}\bigl((0,+\infty)\bigr) -\sE_{T}\bigl((-\infty,0)\bigr).
\end{equation*}
Clearly, the negative and positive subspaces of this involution
coincide with the corresponding spectral subspaces of $T$:
\begin{equation*}
\got H'_{-} = \ran \sE_{T}((-\infty,0)) \quad\text{and}\quad \got
H'_{+} = \ran \sE_{T}((0,+\infty)).
\end{equation*}
Below we will show that in some cases Theorem~\ref{relemma} allows
one to determine the spectral angle of the product $J'J$ where $J$
is another involution on $\fH$. The norm of the difference between
the orthogonal projections onto the corresponding positive (or
negative) subspaces of $J'$ and $J$ is then easily computed by using
\eqref{difPt}.

We study the following two cases.

\begin{hypo}
\label{hypoT} Let $J$ be an involution on the Hilbert space $\fH$.
Assume that $T$ is a self-adjoint operator on $\fH$ such that

$(a)$ $\ker(T)=\{0\}$ and the product $JT$ is accretive

\noindent or

$(b)$ the product $JT$ is strictly accretive.
\end{hypo}

Obviously, if the assumption $(b)$ holds then the assumption $(a)$
holds, too. Therefore, both $(a)$ and $(b)$ assume that
$\ker(T)=\{0\}$. Hence any of these two assumptions implies that the
isometry $J'$ in the polar decomposition~\eqref{polarT} of $T$ is an
involution.

To describe the accretive operators in some more detail we introduce
the following definition.

\begin{defn}
\label{sectorbound} Let $S$ be an accretive operator on the Hilbert
space $\got H$. Then the finite or infinite number
\begin{equation*}
k(S) = \sup_{z\in \cW(S) \backslash\{0\}} \frac{|\im z|}{\re z}
\end{equation*}
is called the \textit{sector bound} of $S$.
\end{defn}

Clearly, if $k(S)$ is finite then $S$ is a sectorial operator (see
\cite[\S V.3.10]{Kato}) with vertex $0$ and semi-angle
$\theta=\arctan k(S)$.

Main result of this section is the following
\begin{thm}
\label{inv} Assume Hypothesis~\ref{hypoT} $(a)$. Let $T = J'|T|$ be
the polar decomposition of~$T$. Then the involutions $J'$ and $J$
are in the acute case, and
\begin{equation}
\label{estspec} \vartheta(U) \leq\frac{1}{2}\arctan k(JT)
\quad\left(\leq\frac{\pi}{4}\right),
\end{equation}
where $U$ is the direct rotation from $J$ to $J'$.
\end{thm}
\begin{proof}
Since $JT$ is accretive and $T=J'|T|$, it follows from
Theorem~\ref{relemma} that the operator $J'J$ is also accretive.
Hence $-1 \notin \cW(J'J)$ and thus by Lemma~\ref{acute} the
involutions $J$ and $J'$ are in the acute case.

If $k(JT)=0$ then $\cW(JT)$ is a subset of the real axis which means
that $J T$ is a symmetric operator. This implies $J\smile T$ since
$T$ is self-adjoint. Hence $J'\smile J$ (see, e.g., \cite[Lemma
VI.2.37]{Kato}) and thus $J=J'$ by Lemma~\ref{remcom}. In this case
estimate~\eqref{estspec} is trivial since $\vartheta(U) = 0$.

Further, assume that $k(JT)>0$. Set
\begin{equation*}
\varphi = \frac{\pi}{2} - \arctan k(JT), \quad \varphi \in
[0,{\pi}/{2}),
\end{equation*}
and observe that the operators $GT$ and $G^*T^*$ with
$G=e^{i\varphi}J$ are both accretive. Then by Theorem \ref{relemma}
one concludes that the products $e^{i\varphi} J'J$ and
$e^{-i\varphi} J'J$ are also accretive operators. Hence $\cW(J'J)$
is a subset of the closed sector
$$
\left\{z\in\bbC\,\left|\,\,\,|\arg z|\leq
\dfrac{\pi}{2}-\varphi\right.\right\}.
$$
Then from the inclusion $\spec(J'J) \subset \overline{\cW(J'J)}$  it
follows that the spectral angle of the unitary operator $J'J$
satisfies
\begin{equation}
\label{varthJJ} \vartheta(J'J) \le \arctan k( JT ).
\end{equation}
Now \eqref{estspec} follows immediately from \eqref{varthJJ} and
\eqref{uhalf}, completing the proof.
\end{proof}

In the two following statements we present some uniqueness results
concerning the involution $J'$ referred to in Theorem \ref{inv}.

\begin{thm}
\label{Juniq} Assume Hypothesis~\ref{hypoT} $(a)$. Let $\widetilde
J'$ be an involution on $\got H$ such that
\begin{equation*}
\mbox{$(i)$ \, $\widetilde J'$ and $J$ are in the acute case, \,
$(ii)$ \, $\widetilde J' \smile T$, \, and \, $(iii)$ \, $\widetilde
J' \ne J'$,}
\end{equation*}
where $J'$ is the involution in the polar decomposition of $T$. Then
\begin{equation}
\label{thprime} \vartheta(\widetilde U) \ge \frac{\pi}{2} -
\frac{1}{2} \arctan k(JT) \quad\left(\geq\frac{\pi}{4}\right),
\end{equation}
where $\widetilde U$ is the direct rotation from $J$ to $\widetilde
J'$.
\end{thm}

\begin{thm}
\label{Juniq1} Assume Hypothesis~\ref{hypoT} $(b)$. Let $T = J' |T|$
be the polar decomposition of~$T$. Then $J'$ is a unique involution
on $\got H$ such that
\begin{equation*}
\mbox{$(i)$ \, $J$ and $J'$ are in the acute case, \, $(ii)$ \,
$J'\smile T$, \, and \, $(iii)$ \, $\vartheta(U) \leq
\dfrac{\pi}{4}$,}
\end{equation*}
where $U$ is the direct rotation from $J$ to $J'$.
\end{thm}

\begin{proof}[Proof of Theorem \ref{Juniq}]
For a proof by contradiction suppose that instead of
\eqref{thprime} the opposite inequality holds. Then by
\eqref{uhalf} in Remark \ref{star} we have
\begin{equation}
\label{thp1} \vartheta(\widetilde J' J) < \pi - \arctan k(JT).
\end{equation}
Similarly, Theorem \ref{inv} yields
\begin{equation}
\label{thp2} \vartheta(J J') \le \arctan k(J T).
\end{equation}
By \eqref{thp1} and \eqref{thp2} Lemma~\ref{UU} implies that
\begin{equation*}
\vartheta(\widetilde J' J') = \vartheta((\widetilde J' J)  (J J'))
\le \vartheta(\widetilde J' J) + \vartheta(J J') < \pi.
\end{equation*}
In particular, this means that  $-1\notin \spec(\widetilde{J'} J')$
which proves that the involutions $J'$ and $\widetilde{J'}$ are in
the acute case.

By hypothesis $\widetilde J'$ commutes with $T$ and $J'$ is the
isometry in the polar decomposition of $T$. Hence \cite[Lemma
VI.2.37]{Kato} implies $ \widetilde J' \smile J'$. Then from Lemma
\ref{remcom} it follows that $\widetilde J' = J'$ which contradicts
the assumption $(iii)$. Therefore $\vartheta(\widetilde{U})$
satisfies \eqref{thprime} completing the proof.
\end{proof}

\begin{proof}[Proof of Theorem \ref{Juniq1}]

Arguing by contradiction, suppose that there is an involution
$\widetilde J'$ distinct from $J'$ and such that conditions
$(i)$--$(iii)$ are satisfied. In particular, this implies  that
$\vartheta(\widetilde J' J)\le \pi/2$ and hence
\begin{equation*}
\re \lal x, J \widetilde J'  x \ral \ge 0 \quad\text{for all}\quad
x\in \got H.
\end{equation*}
Since $J T$ is strictly accretive and $T = J' |T|$, by
Theorem~\ref{relemma} the operator $JJ'$ is also strictly accretive,
that is,
\begin{equation*}
\re \lal x, J J' x \ral > 0 \quad\text{for all}\quad x\in \got H,
\quad x\ne 0.
\end{equation*}
Therefore,
\begin{equation}
\label{strineq0} \re \lal x, J J' x \ral + \re \lal x, J \widetilde
J'  x \ral > 0 \text{\, for all \,} x\in\got H, \quad x\ne 0.
\end{equation}

Now assume that there is $y\in \ker(I + \widetilde J' J')$ such that
$y\ne 0$. Then applying $\widetilde J'$ to both parts of the
equality $y + \widetilde J' J' y = 0$ yields $J' y + \widetilde J' y
= 0$. Hence
\begin{equation*}
\re \lal y, J J' y \ral + \re \lal y, J \widetilde J' y \ral = 0,
\end{equation*}
and it follows from \eqref{strineq0} that $y = 0$. This proves that
$\ker(I + \widetilde J' J') = \{0\}$, i.e. the involutions $
\widetilde J'$ and  $J'$ are in the acute case.

Clearly, $ \widetilde J' \smile J'$ since by  hypothesis $\widetilde
J'$ commutes with $T$ and $J'$ is the isometry in the polar
decomposition of $T$ (see \cite[Lemma VI.2.37]{Kato}). Hence, by
Lemma~\ref{remcom} $\widetilde J' = J'$ which contradicts the
assumption that $\widetilde J'$ is distinct from  $J'$.

The proof is complete.
\end{proof}

\section{An extension of the Davis-Kahan $\boldsymbol{\tan2\Theta}$ Theorem.
\\ Proof of Theorem 1} \label{Stan2t}
%% Proof of Theorem \ref{Th1}

Throughout this section we adopt the following hypothesis.

\begin{hypo}
\label{Hyp2t} Given a self-adjoint operator $A$ on the Hilbert space
$\got H$ assume that
\begin{equation}
\label{Amu} \ker ( A - \mu ) =\{ 0 \} \quad\text{for some}\quad
\mu\in\mathbb{R}.
\end{equation}
Let $V$ be a symmetric operator on $\got H$ such that
\smallskip

\begin{itemize}
\item[$(i)$] $\dom(A)\subset\dom(V)$,

\item[$(ii)$] $V \frown J$ where
$J = \sE_A((\mu,+\infty)) - \sE_A((-\infty,\mu))$,
\end{itemize}
and
\begin{itemize}
\item[$(iii)$] the closure $L=\overline{L_0}$ of the operator $L_0 = A +
V$ with $\dom(L_0)=\dom(A)$ is a self-adjoint operator.
\end{itemize}
\end{hypo}

Under this hypothesis the product $J(L-\mu)$ appears to be a
strictly accretive operator. Moreover, the sector bound
$k\bigl(J(L-\mu)\bigr)$ admits an explicit description in terms of
the perturbation~$V$.

\begin{lem}
\label{remHyp2t} Assume Hypothesis \ref{Hyp2t}. Then $J(L-\mu)$ is a
strictly accretive operator and
\begin{equation}
\label{kJLm} k\bigl(J(L-\mu)\bigr)= \sup\limits_{\mbox{\scriptsize
$\begin{array}{cc} x\in \dom(A)\\ \|x\|=1\end{array}$}} \frac{|\lal
x,JVx\ral|}{\lal x,|A-\mu|x\ral}.
\end{equation}
\end{lem}
\begin{proof}
Obviously, under Hypothesis \ref{Hyp2t}
\begin{equation*}
J(A-\mu)=|A - \mu| > 0.
\end{equation*}
Hence by items $(ii)$ and $(iii)$ of this hypothesis
\begin{equation}
\label{r1} \re\lal x, J(A+V-\mu)x\ral =\lal x,|A-\mu| x\ral \text{\,
for all \,}x\in \dom(A).
\end{equation}
Pick up an arbitrary $y\in \dom(L)$. By the assumption
\textit{(iii)} it follows that there exists a sequence of vectors
$y_n\in\dom(A)$ such that $y_n\to y$ and $L_0 y_n \to L y$ as
$n\to\infty$, and thus
\begin{equation}
\label{lr} \re \lal y_n, J(L_0-\mu) y_n\ral \to \re\lal y,
J(L-\mu)y\ral \text{\, as \,} n\to \infty.
\end{equation}
Then \eqref{r1} and \eqref{lr} imply $\re\lal y, J(L-\mu) y\ral\ge
0$. Moreover, $y\in \ker(|A-\mu|)\subset\dom(A)$ whenever $\re\lal
y, J(L-\mu) y\ral= 0$. Taking into account that $\ker (|A-\mu|)=\ker
(A-\mu) = \{0\}$, one infers that
\begin{equation*}
\re\lal y, J(L-\mu) y\ral > 0 \text{\, for all non-zero \,}
y\in{\dom(L)},
\end{equation*}
which means that the operator $J(L-\mu)$ is strictly accretive.

Now observe that
\begin{equation}
\label{kbol} k\bigl(J(L-\mu)\bigr)\geq\varkappa,
\end{equation}
where
\begin{equation}
\label{kappa} \varkappa= \sup\limits_{\mbox{\scriptsize
$\begin{array}{cc} x\in \dom(A)\\ \|x\|=1\end{array}$}} \frac{|\lal
x,JVx\ral|}{\lal x,|A-\mu|x\ral}.
\end{equation}
Indeed,
\begin{align*}
k\bigl(J(L-\mu)\bigr)&= \sup\limits_{\mbox{\scriptsize$
\begin{array}{cc} x\in \dom(L)\\ \|x\|=1\end{array}$}}
\frac{|\im\lal x,J(L-\mu)x\ral|}{\re\lal x,J(L-\mu)x\ral}\\
&\geq \sup\limits_{\mbox{\scriptsize$
\begin{array}{cc} x\in \dom(A)\\ \|x\|=1\end{array}$}}
\frac{|\im\lal x,J(A+V-\mu)x\ral|}{\re\lal x,J(A+V-\mu)x\ral}
\end{align*}
since by Hypothesis \ref{Hyp2t} $(iii)$ $
\dom(A)\subset\dom(L)\text{\, and \,}L|_{\dom(A)}=A+V. $ Then
\eqref{kbol} holds by \eqref{r1} since Hypothesis \ref{Hyp2t} $(ii)$
implies
\begin{equation}
\label{ImL} \im\lal x,J(A+V-\mu)x\ral=\lal x,JVx\ral \text{\, for
any \,}x\in\dom(A).
\end{equation}

Clearly,  if $\varkappa=\infty$ then \eqref{kJLm} follows
immediately from inequality \eqref{kbol}. If $\varkappa$ is finite,
then by \eqref{r1} and \eqref{ImL} from \eqref{kappa} we have
\begin{equation*}
|\im\lal x,J(L_0-\mu)x\ral|\leq\varkappa\re\lal x,J(L_0-\mu)x\ral
\text{\, for any \,}x\in\dom(L_0)=\dom(A).
\end{equation*}
Since $L$ is the closure of $L_0$, by continuity of the inner
product the same inequality holds for $L$, that is,
\begin{equation*}
|\im\lal x,J(L-\mu)x\ral|\leq\varkappa\re\lal x,J(L-\mu)x\ral
\text{\, for any \,}x\in\dom(L).
\end{equation*}
In particular, this means that
\begin{equation}
\label{kmen} \sup\limits_{\mbox{\scriptsize$
\begin{array}{cc} x\in \dom(L)\\ \|x\|=1\end{array}$}}
\frac{|\im\lal x,J(L-\mu)x\ral|}{\re\lal x,J(L-\mu)x\ral}
=k\bigl(J(L-\mu)\bigr)\leq \varkappa.
\end{equation}
Now combining \eqref{kbol}, \eqref{kappa}, and \eqref{kmen}
completes the proof.
\end{proof}

\begin{rem}
\label{reminv} Since $J(L-\mu)$ is a strictly accretive operator,
the isometry $J'$ in the polar decomposition $L-\mu = J'|L-\mu|$ is
an involution. Clearly, it reads
\begin{equation*}
%\label{jprimed}
J' = \sE_{L}((\mu,+\infty)) - \sE_{L}((-\infty,\mu)).
\end{equation*}
\end{rem}

\begin{thm}
\label{dkn} Assume Hypothesis \ref{Hyp2t}. Let $L-\mu = J'|L-\mu|$
be the polar decomposition of $L-\mu$. Then the involutions $J$ and
$J'$ are in the acute case, and
\begin{equation}
\label{tg2} \vartheta(U)\le \frac{1}{2}\arctan\Bigl(
\sup\limits_{\mbox{\scriptsize$\begin{array}{cc} x\in \dom(A)\\
\|x\|=1\end{array}$}} \frac{|\lal x,JVx\ral|}{\lal x,|A-\mu| x\ral}
\Bigr) \quad\left(\leq\frac{\pi}{4}\right),
\end{equation}
where $U$ is the direct rotation from $J$ to $J'$. Moreover, $J'$ is
a unique involution on $\fH$ with the properties
\begin{equation}
\label{upart} \mbox{$(i)$ \, $J'$ and $J$ are in the acute case, \,
$(ii)$ \, $J' \smile L$, \, and \, $(iii)$ \, $\vartheta(U) \le
\dfrac{\pi}{4}$.}
\end{equation}
The spectral angle of the direct rotation $\widetilde{U}$ from $J$
to any other involution $\widetilde{J'}$ distinct from $J'$ and
satisfying $(i)$ and $(ii)$ is bounded from below as follows
\begin{equation}
\label{thprime1} \vartheta(\widetilde{U}) \ge \frac{\pi}{2} -
\frac{1}{2} \arctan\Bigl(
\sup\limits_{\mbox{\scriptsize$\begin{array}{cc} x\in \dom(A)\\
\|x\|=1\end{array}$}} \frac{|\lal x,JVx\ral|}{\lal x,|A-\mu| x\ral}
\Bigr) \quad\left(\geq\frac{\pi}{4}\right).
\end{equation}
\end{thm}
\begin{proof}
The operators $J$ and $T = L-\mu$ satisfy Hypothesis~\ref{hypoT} (b)
(and hence, Hypothesis~\ref{hypoT} (a)). Then the assertion is
proven simply by combining Theorems \ref{inv}, \ref{Juniq}, and
\ref{Juniq1} with Lemma~\ref{remHyp2t}.
\end{proof}

With Theorem \ref{dkn} one can easily prove Theorem \ref{Th1}.

\begin{proof}[Proof of Theorem \ref{Th1}]
Pick up arbitrary $\mu, \nu \in(\sup\sigma_-,\inf\sigma_+)$, $\mu <
\nu$. Clearly, Hypothesis~\ref{Hyp2t} holds for both $\mu$ and $\nu$
with the same involution $J = \sE_{A}(\sigma_{+}) -
\sE_{A}(\sigma_{-})$. By Remark \ref{reminv} the isometries
$J'_{\mu}$ and $J'_{\nu}$ in the polar decompositions $L-\mu =
J'_{\mu}|L-\mu|$  and $L-\nu = J'_{\nu}|L-\nu|$ are involutions. By
Theorem \ref{dkn} the involutions $J$ and $J'_{\mu}$ are in the
acute case, $J'_{\mu} \smile L$, and $\vartheta(U_{\mu}) \le \pi/4$
where $U_{\mu}$ is the direct rotation from $J$ to $J'_{\mu}$. The
same holds for $J'_{\nu}$ and the corresponding direct rotation
$U_{\nu}$ from $J$ to $J'_{\nu}$. Therefore, \eqref{upart} is
satisfied for both $J' = J'_{\mu}$ and $J' = J'_{\nu}$. Hence,
Theorem \ref{dkn} implies $J'_{\mu} = J'_{\nu}$ which by Remark
\ref{reminv} yields $\sE_L\bigl((\mu,\nu)\bigr)=0$. Since $\mu, \nu
\in(\sup\sigma_-,\inf\sigma_+)$ are arbitrary, one then concludes
that $\sE_L\bigl( (\sup \sigma_{-}, \inf \sigma_{+} ) \bigr) =0$,
and thus the interval $(\sup\sigma_-,\inf\sigma_+)$ belongs to the
resolvent set of $L$. Hence,
\begin{equation*}
J'_{\mu} = \sE_{L}(\sigma'_{+}) - \sE_{L}(\sigma'_{-})
\quad\text{for all}\quad \mu \in(\sup\sigma_-,\inf\sigma_+),
\end{equation*}
where $\sigma'_{-}$ and $\sigma'_{+}$ are the parts of the spectrum
of $L$ in the intervals $(-\infty, \sup \sigma_{-}]$ and $[\inf
\sigma_{+} , +\infty)$, respectively. Since $J'_{\mu}$ does not
depend on $\mu \in(\sup\sigma_-,\inf\sigma_+)$, the direct rotation
$U_{\mu}$ does not, too. Then estimate~\eqref{tg2} of
Theorem~\ref{dkn} yields
\begin{equation}
\label{infest} \tan 2 \vartheta(U) \le \inf_{\sup {\sigma_{-}} < \mu
< \inf {\sigma_{+}} }
\sup\limits_{\mbox{\scriptsize$\begin{array}{cc} x\in \dom(A)\\
\|x\|=1\end{array}$}} \frac{|\lal x,JVx\ral|}{\lal x,|A-\mu| x\ral},
\end{equation}
where $U$ is the direct rotation from the involution
$\sE_{A}(\sigma_{+}) - \sE_{A}(\sigma_{-})$  to the involution
$\sE_{L}(\sigma'_{+}) - \sE_{L}(\sigma'_{-})$. Now inequality
\eqref{infest} proves the bound~\eqref{Estin} by taking into account
\eqref{difPt} in Remark \ref{star}. The proof is complete.
\end{proof}

\begin{ex}
\label{ksharp} Let $\cD_{a} = \mathbb{R}\backslash(-a,a)$ for some
$a\ge 0$. Given $\varkappa\geq 0$ assume that $A$ and $V$ are
operators on the Hilbert space $\got H = L^2(\cD_a)$ defined by
\begin{equation}
\label{AVEX}
\begin{aligned}
& \quad (A x)(t) = |t| x(-t), \quad (V x)(t) = \varkappa\, t\, x(t),
\quad
t\in\cD_a,\\
& \dom(A)=\dom(V) =\bigl\{ x \in \got H \,|\,\, \int_{\cD_a} t^2
|x(t)|^2 d t < +\infty\bigr\}.
\end{aligned}
\end{equation}
Both $A$ and $L=A+V$ are self-adjoint operators. The spectrum of the
operator $A$ is purely absolutely continuous. For $a>0$ it consists
of the two disjoint components $\sigma_{-} = (-\infty,-a]$ and
$\sigma_{+} = [a,+\infty)$ and for $a=0$ it covers the whole real
axis. Obviously, the isometry $J$ in the polar decomposition $A =
J|A|$ is the parity operator, $(Jx)(t) = x(-t)$, $x\in\fH$, and the
absolute value of $A$ is given by $(|A|x)(t) = |t|x(t)$,
$x\in\Dom(A)$. Clearly, $J$ is an involution on $\got H$ such that
$J\smile A$ and $J\frown V$. Therefore, for $a > 0$ (resp. for
$a=0$) the operators $A$ and $V$ satisfy the hypothesis of
Theorem~\ref{Th1} (resp. the hypothesis of Theorem~\ref{dkn} for
$\mu=0$).

Our analysis of the subspace perturbation problem involving $A$ and
$V$ given by \eqref{AVEX} is divided into three parts below.

$(i)$ For any $x\in \dom(A)$, $\|x\|=1$, we have
\begin{align}
|\lal x, JV x\ral| &  =  \Bigl| \int_{\cD_a} \varkappa t
\overline{x(t)} x(-t) d t \Bigr| \notag
 \\  & \le  \varkappa \int_{\cD_a} |t| |\overline{x(t)} x(-t)| d t \notag
 \\ & \le  \varkappa \int_{\cD_a} |t| \frac{|x(-t)|^2
+ |x(t)|^2}{2}d t \label{zz2}
 \\  & = \varkappa \int_{\cD_a} |t| |x(t)|^2 d t \notag
 \\  & = \varkappa \lal x,|A| x\ral,\notag
\end{align}
Moreover, if $x\in \dom(A)$ is such that $x(-t) = i
\mathop{\mathrm{sign}}(t) x(t)$ then inequalities in \eqref{zz2}
turn into equalities. Hence, by taking this into account,
\eqref{zz2} implies
\begin{equation}
\label{kapkap} \sup\limits_{\mbox{\scriptsize$\begin{array}{cc} x\in
\dom(A)\\ \|x\|=1\end{array}$}} \frac{|\lal x,JVx\ral|}{\lal x,|A|
x\ral} =\varkappa.
\end{equation}
An explicit evaluation of the involution $J'=\sE_L((+\infty,0)) -
\sE_L((-\infty,0))$ by using the polar decomposition $L=J'|L|$
yields
\begin{equation}
\label{JJxa} (J'J x)(t) = \frac{1}{\sqrt{1+\varkappa^2}} x(t) +
\mathop{\mathrm{sign}}(t) \frac{\varkappa}{\sqrt{1+\varkappa^2}}
x(-t).
\end{equation}
From \eqref{JJxa} it follows by inspection that the spectrum of the
unitary operator $J'J$ consists of the two mutually conjugate
eigenvalues,
\begin{equation*}
\spec(J'J) = \left\{ \frac{1-i\varkappa}{\sqrt{1+\varkappa^2}},
\frac{1+i\varkappa}{\sqrt{1+\varkappa^2}} \right\}.
\end{equation*}
This implies that $\vartheta(J'J) = \arctan \varkappa$ and then the
spectral angle of the direct rotation $U$ from $J$ to $J'$ is equal
to $\vartheta(U)=\dfrac{1}{2}\arctan \varkappa$. Combining this with
\eqref{kapkap} yields that in the case under consideration
\begin{equation}
\label{tg11} \vartheta(U)=\frac{1}{2}\arctan\Bigl(
\sup\limits_{\mbox{\scriptsize$\begin{array}{cc} x\in \dom(A)\\
\|x\|=1\end{array}$}} \frac{|\lal x,JVx\ral|}{\lal x,|A| x\ral}
\Bigr)\quad\text{for any}\quad a\geq 0.
\end{equation}

$(ii)$ Now set $\widetilde{J'}= -J'$. Clearly,
$\vartheta(\widetilde{J'} J) = \pi - \vartheta(J'J)$ and thus the
spectral angle of the direct rotation $\widetilde{U}$ from $J$ to
$\widetilde{J'}$ reads
\begin{equation}
\label{zz3} \vartheta(\widetilde{U}) = \frac{\pi}{2} - \frac{1}{2}
\arctan\Bigl( \sup\limits_{\mbox{\scriptsize$\begin{array}{cc} x\in
\dom(A)\\ \|x\|=1\end{array}$}} \frac{|\lal x,JVx\ral|}{\lal x,|A|
x\ral} \Bigr).
\end{equation}
Notice that the involution $\widetilde{J'}$ commutes with $L$ since
$J'$ does. By \eqref{JJxa} it follows that $\ker(I-J'J)=\{0\}$
whenever $\varkappa\neq 0$. Hence $\ker(I +
\widetilde{J'}\,J)=\{0\}$ whenever $\varkappa\neq0$ which means that
for $\varkappa>0$ the involutions $J$ and $\widetilde{J'}$ are in
the acute case.

$(iii)$ For $a>0$ we have
\begin{equation}
\label{zz4} \inf_{|\mu| < a}
\sup\limits_{\mbox{\scriptsize$\begin{array}{cc} x\in \dom(A)\\
\|x\|=1\end{array}$}} \frac{|\lal x,J V x\ral|}{\lal
x,|A-\mu|x\ral}\leq \sup\limits_{\mbox{\scriptsize$\begin{array}{cc}
x\in \dom(A)\\ \|x\|= 1\end{array}$}} \frac{|\lal x,JVx\ral|}{\lal
x,|A| x\ral}.
\end{equation}
Since $\sin\vartheta(U)=\left\|\sE_L\bigl((-\infty,-a]\bigr) -
\sE_A\bigl((-\infty,-a]\bigr)\right\|$, by Theorem \ref{Th1} the
strict inequality in \eqref{zz4} implies
\begin{equation*}
\vartheta(U)<\dfrac{1}{2}\arctan\Bigl(
\sup\limits_{\mbox{\scriptsize $\begin{array}{cc} x\in \dom(A)\\
\|x\|= 1\end{array}$}} \frac{|\lal x,JVx\ral|}{\lal x,|A|
x\ral}\Bigr),
\end{equation*}
which contradicts \eqref{tg11}. Hence only the equality sign in
\eqref{zz4} is allowed and thus
\begin{equation}
\label{zz5} \inf_{|\mu| < a} \sup\limits_{\mbox{\scriptsize
$\begin{array}{cc} x\in \dom(A)\\ \|x\|=1\end{array}$}} \frac{|\lal
x,J V x\ral|}{\lal x,|A-\mu|x\ral}=\sup\limits_{\mbox{\scriptsize
$\begin{array}{cc} x\in \dom(A)\\ \|x\|=1\end{array}$}} \frac{|\lal
x,J V x\ral|}{\lal x,|A|x\ral}.
\end{equation}

\end{ex}

\begin{rem}
\label{remTh1} Example \ref{ksharp} shows the following.
\begin{itemize}
\item[$(i)$] Estimate \eqref{tg2}
of Theorem \ref{dkn} is sharp. This is proven by equality
\eqref{tg11}.
\item[$(ii)$] Estimate \eqref{thprime1}
of the same theorem is sharp. This is proven by equality
\eqref{zz3}.
\item[$(iii)$] Estimate \eqref{Estin} of Theorem \ref{Th1} is sharp.
This is proven by combining equalities \eqref{tg11} and \eqref{zz5}.
\end{itemize}
\end{rem}

The celebrated sharp estimate for the operator angle between the
spectral subspaces $\ran \sE_A(\sigma_{-})$ and $\ran
\sE_L(\sigma'_{-})$ known as the Davis-Kahan $\tan2\Theta$ Theorem
\cite{DK} (cf. \cite{KMM5}) appears to be a simple corollary to
Theorem \ref{Th1}.

\begin{thm}[The Davis--Kahan $\boldsymbol{\tan2\Theta}$ Theorem]
\label{DK2T} Given a self-adjoint operator $A$ on the Hilbert space
$\got H$ assume that
\begin{equation*}
\spec(A) = \sigma_{-}\cup \sigma_{+}, \quad d = \dist(\sigma_{-},
\sigma_{+}) > 0, \quad \text{and} \quad \sup\sigma_{-} < \inf
\sigma_{+}.
\end{equation*}
Suppose that a bounded self-adjoint operator $V$ on $\got H$ is
off-diagonal with respect to the decomposition $\got H = \ran
E_{A}(\sigma_{-})\oplus \ran E_{A}(\sigma_{+})$. Then the spectrum
of $L = A + V$ consists of two disjoint components $\sigma'_{-}$ and
$\sigma'_{+}$ such that
\begin{equation*}
\sigma'_{-} \subset (-\infty, \sup{\sigma_{-}}] \quad
\text{and}\quad \sigma'_{+} \subset [\inf \sigma_{+}, +\infty),
\end{equation*}
and
\begin{equation}
\label{EstinDK} \| \sE_{L}(\sigma_{-}') - \sE_A(\sigma_{-})\| \le
\sin\Bigl(\frac{1}{2} \arctan \frac{2 \| V \|}{d}\Bigr).
\end{equation}
\end{thm}

\begin{proof}
Hypothesis of Theorem \ref{Th1} is satisfied and thus we only need
to prove the estimate \eqref{EstinDK}. Set $\mu_0 =
\dfrac{1}{2}(\sup\sigma_{+} + \inf\sigma_{-})$. Clearly,
\begin{align*}
\nonumber \inf_{\sup {\sigma_{-}} < \mu < \inf {\sigma_{+}} }
\sup\limits_{\mbox{\scriptsize$\begin{array}{cc} x\in
\dom(A)\\\nonumber \|x\|=1\end{array}$}} \frac{|\lal
x,JVx\ral|}{\lal x,|A-\mu| x\ral} &\le \sup
\limits_{\mbox{\scriptsize$\begin{array}{cc} x\in \dom(A)\\
\|x\|=1\end{array}$}} \frac{|\lal x,JVx\ral|}{\lal x,|A-\mu_0|
x\ral} \\\nonumber &\le \sup
\limits_{\mbox{\scriptsize$\begin{array}{cc} x\in \dom(A)\\
\|x\|=1\end{array}$}} \frac{\|V\|}{\lal
x,|A-\mu_0| x\ral} \\
%\label{Vd2}
& \le \frac{2\|V\|}{d},
\end{align*}
which immediately implies \eqref{EstinDK} by taking into account
\eqref{Estin}.
\end{proof}

\section{Proof of Theorem 2}
\label{S3isl}
%%% Proof of Theorem \ref{main0}

In the proof of the main result of this section we will use some
auxiliary statements. We start with the following lemma.

\begin{lem}
\label{wt} Let $T$ be a densely defined operator on a Hilbert space
$\fH$ with $\dim(\fH)\geq n$ for some $n\in\bbN$. Assume that
$\mathfrak{t}(x,y)$ is a sesquilinear form on $\fH$ such that
\begin{equation*}
%\label{Tincl}
\dom(T)\subset\dom(\mathfrak{t}) \quad\text{and}\quad
\mathfrak{t}(x,y)=\lal x,Ty\ral \quad\text{for any}\quad
x,y\in\dom(T).
\end{equation*}
Suppose that there are orthogonal projections $P_i\neq0$,
$i=1,2,\ldots,n,$ on $\fH$ with the properties
\begin{equation*}
%\label{Pinc} &
P_i P_j = 0 \text{\, if \,} i\ne j,  \quad \sum_{i=1}^n P_i = I,
\text{\, and \,} P_i x\in\dom(\mathfrak{t}) \text{\, whenever
\,}x\in\dom(\mathfrak{t}).
\end{equation*}
Let $\mathscr{E}$ be a set of ordered $n$-element orthonormal
systems in $\got H$ defined by
\begin{equation*}
\mathscr{E} = \big\{\{e_i\}_{i=1}^n \subset
\dom(\mathfrak{t})\,\big|\,\, e_i\in\ran P_i \text{\, and \,}
\|e_i\| = 1 \text{\, for all \,} i=1,2,...n\big\}.
\end{equation*}
Then
\begin{equation}
\label{inc} \cW(T) \subset \bigcup_{\mathbf{e}\in \mathscr{E}}
\cW(\got t^{\mathbf{e}}),
\end{equation}
where for any $\mathbf{e}\in\mathscr{E}$ the $n\times n$ matrix
$\got t^{\mathbf{e}}$ is given by
\begin{equation*}
(\got t^{\mathbf{e}})_{ij}= \mathfrak{t}(e_i,e_j) \text{\, with \,}
e_i,e_j\in\mathbf{e}, \quad i,j=1,2,\ldots,n.
\end{equation*}
If, in addition, $\dom(\mathfrak{t})=\dom(T)$ then
\begin{equation}
\label{inc1} \cW(T) = \bigcup_{\mathbf{e}\in \mathscr{E}} \cW(\got
t^{\mathbf{e}}).
\end{equation}
\end{lem}

\begin{proof}
By hypothesis $\overline{\dom(T)}=\fH$ and hence
$\overline{\dom(\mathfrak{t})}=\fH$, too. Therefore there exists
$y\in\dom(\mathfrak{t})$ such that $P_iy\neq0$ for all
$i=1,2,\ldots,n$. Set $e_i=\dfrac{P_iy}{\|P_iy\|}$. Taking into
account that by hypothesis $P_iy\in\dom(\mathfrak{t})$ and thus
$e_i\in\dom(\mathfrak{t})$, $i=1,2,\ldots,n$, one concludes that
$\{e_{i}\}_{i=1}^n\in \mathscr{E}$. Hence, the set $\mathscr{E}$ is
non-empty.

Assume that $z\in \cW(T)$. Then there exists $x\in\dom(T)$ such that
$\lal x, T x\ral =z$ and $\|x\|=1$. Pick up an arbitrary
$\mathbf{f}=\{f_i\}_{i=1}^n\in\mathscr{E}$ and define the
orthonormal system $\mathbf{g} = \{ g_i \}_{i=1}^n$ by
\begin{equation*}
  g_i = \left\{
  \begin{aligned}
    \frac{P_i x}{\|P_i x\|}, \quad & \|P_i x\| \ne 0, \\
    f_i, \quad & \|P_i x\| = 0.
    \end{aligned}
  \right.
\end{equation*}
Obviously, $\mathbf{g}\in \mathscr{E}$ and
\begin{equation*}
\sum_{i,j=1}^n \mathfrak{t}(g_i,g_j) \|P_i x\| \| P_j x\|=\lal x,
Tx\ral=z,
\end{equation*}
which implies $z\in \cW(\got t^{\mathbf{g}})$ since
$\displaystyle\sum\nolimits_{i=1}^n \|P_i x\|^2 = \| x\|^2 = 1$.
This proves the inclusion \eqref{inc}.

To prove the converse inclusion in the case where
$\dom(\mathfrak{t})=\dom(T)$, pick up an arbitrary
$\mathbf{h}=\{h_i\}_{i=1}^n\in \mathscr{E}$ and assume that $z\in
\cW(\got t^{\mathbf{h}})$. Then there are $\alpha_i\in \mathbb{C}$,
$i = 1,2,...,n,$ such that
\begin{equation*}
z= \sum_{i,j=1}^n \mathfrak{t}(h_i,h_j)\alpha_i \overline{\alpha}_j,
\quad \sum_{i=1}^n|\alpha_i|^2 = 1.
\end{equation*}
Set $x=\displaystyle\sum\nolimits_{i=1}^n \alpha_i h_i$. Clearly,
$\|x\| = 1$ and $x\in\dom(\mathfrak{t})=\dom(T)$. Hence
$z=\mathfrak{t}(x,x)=\lal Tx,x\ral$. This yields $z\in \cW(T)$ and
hence $\cW(\got t^{\mathbf{h}})\subset \cW(T)$. One then concludes
that
\begin{equation*}
\bigcup_{\mathbf{e}\in \mathscr{E} } \cW(\got t^{\mathbf{e}})\subset
\cW(T)
\end{equation*}
and hence~\eqref{inc1} holds, completing the proof.
\end{proof}

The next simple result on the numerical range of a $2\times2$
numerical matrix is well known (see, e.g., \cite[Lemma 1.1-1]{Rao}).

\begin{lem}
\label{ellipse} Given numbers $\alpha > 0$, $\beta > 0$, and $\gamma
\in \mathbb{C}$ let $M$ be a $2\times 2$ matrix of the form
\begin{equation*}
  M =\left(
  \begin{array}{rr}
    \alpha & -\overline{\gamma} \\
    \gamma & \beta
  \end{array}\right).
\end{equation*}
The matrix $M$ is strictly accretive and its sector bound reads
\begin{equation*}
k(M)=\frac{|\gamma|}{\sqrt{\alpha \beta}}.
\end{equation*}
The numerical range $\cW(M)$ is a (possibly degenerate) elliptical
disc with foci at the eigenvalues of $M$.
\end{lem}

Now we are in a position to prove the main statement of the section.
We only recall that by a finite gap of a closed set $\sigma\subset\bbR$ one
understands an \emph{open} finite interval on the real axis that
does not intersect this set but both its ends belong to $\sigma$.

\begin{thm}
\label{trio1} Given a self-adjoint operator $A$ on the Hilbert space
$\got H$ assume that
\begin{equation*}
\spec(A) = \sigma_{-}\cup\sigma_{+},  \quad
\dist(\sigma_{+},\sigma_{-}) = d > 0, \quad  \text{and} \quad
\sigma_{-}\subset\Delta,
\end{equation*}
where $\Delta=(\alpha,\beta)$, $\alpha<\beta$, is a finite gap of
$\sigma_+$. Assume in addition that $V$ is a bounded self-adjoint
operator on $\got H$ anticommuting with
$J=\sE_A(\sigma_+)-\sE_A(\sigma_-)$ and such that
\begin{equation}
\label{Vdd} \| V \| < \sqrt{d(|\Delta| - d) },
\end{equation}
where $|\Delta|=\beta-\alpha$ denotes the length of the gap
$\Delta$. Then the spectrum of $L= A+V$ consists of two disjoint
components $\sigma'_{-}$ and $\sigma'_{+}$ such that inclusions
\eqref{incl2} hold with $\delta_\mp$ given by \eqref{incl2a},
\eqref{incl2b} and the involutions $J$ and $J'=\sE_{L}(\sigma'_{+})
- \sE_{L}(\sigma'_{-})$ are in the acute case. The spectral angle of
the direct rotation $U$ from $J$ to $J'$ satisfies the bound
\begin{equation}
\label{thebest} \vartheta(U) \leq\frac{1}{2}\arctan
\kappa\bigl(\|V\|\bigr) \quad\Bigl(\leq \frac{\pi}{4}\Bigr),
\end{equation}
where the function $\kappa(v)$ is defined for $0\leq
v<\sqrt{d(|\Delta|-d)}$ by
\begin{equation}
\label{kappaV} \kappa(v)=\left\{\begin{array}{cl} \displaystyle
\frac{2v}{d} & \text{ if \,} v \le
\displaystyle\sqrt{\frac{d}{2}\left(\frac{|\Delta|}{2}-d\right)},\\
\displaystyle \frac{v\dfrac{|\Delta|}{2} + \sqrt{d(|\Delta| - d
)\Bigl[ \Bigl(\dfrac{|\Delta|}{2} -  d\Bigr)^2 + v^2 \Bigr]}}
{d(|\Delta|-d) - v^2} & \text{ if \,} v
>\displaystyle\sqrt{\frac{d}{2}\left(\frac{|\Delta|}{2}-d\right)}.
\end{array}\right.
\end{equation}
Moreover, $J'$ is a unique involution on $\fH$ with the properties
\begin{equation*}
\mbox{$(i)$ \, $J'$ and $J$ are in the acute case, \, $(ii)$ \, $J'
\smile L$, \, and \, $(iii)$ \, $\vartheta(U) \le \dfrac{\pi}{4}$.}
\end{equation*}
The spectral angle of the direct rotation $\widetilde{U}$ from $J$
to any involution $\widetilde{J'}$ distinct from $J'$ and satisfying
$(i)$ and $(ii)$ is bounded from below as follows
\begin{equation}
\label{thprime2} \vartheta(\widetilde{U}) \ge \frac{\pi}{2} -
\frac{1}{2} \arctan\kappa(\|V\|).
\end{equation}
\end{thm}

\begin{proof}
Recall that inclusions \eqref{incl2} with $\delta_\mp$ given by
\eqref{incl2a}, \eqref{incl2b} follow from \cite[Theorems 1 (i) and
3.2]{KMM3}. In the proof of the remaining statements one may assume
without loss of generality that the gap $\Delta$ is centered at the
point zero. Under this assumption we set
\begin{equation}
\label{abet} \alpha=-b \text{\, and \,} \beta=b \text{\, with \,}
b=\frac{|\Delta|}{2}.
\end{equation}
Then
\begin{equation}
\label{intervals1} \sigma_+\subset\bbR\setminus(-b,b)\text{\, and
\,}  \sigma_-\subset [-a,a],
\end{equation}
where
\begin{equation*}
%\label{ab}
a=\frac{|\Delta|}{2} - d, \quad 0\le a < b.
\end{equation*}

For $\alpha,\beta$ given by \eqref{abet}, inclusions \eqref{incl2}
imply that the intervals $(-b,-a')$ and $(a',b)$ with
\begin{equation*}
a'=a+\|V\|\tan\left(\frac{1}{2}\arctan\frac{2\|V\|}{a+b}\right)<b
\end{equation*}
are in the resolvent set of $L$. Hence the interval $(a'^2,b^2)$
lies in the resolvent set of $L^2$. Taking into account \eqref{Vdd}
one verifies by inspection that $a'^2\leq a^2 + \|V\|^2<b^2$.
Therefore, the interval \mbox{$(a^2 + \|V\|^2, b^2)$} belongs to the
resolvent set of $L^2$ and the spectral projections
$\sE_{L^2-\mu}\bigl((-\infty,0)\bigr)$ and
$\sE_{L^2-\mu}\bigl((0,\infty)\bigr)$ do not depend on
\begin{equation}
\label{much} \mu\in(a^2 + \|V\|^2, b^2).
\end{equation}
Moreover,
\begin{equation*}
\sE_{L^2-\mu}\bigl((-\infty,0)\bigr)=\sE_L(\sigma'_-), \quad
\sE_{L^2-\mu}\bigl((0,\infty)\bigr) =\sE_L(\sigma'_+),
\end{equation*}
and hence
\begin{equation}
\label{JLm} \sE_{L^2-\mu}\bigl((0,\infty)\bigr)-
\sE_{L^2-\mu}\bigl((-\infty,0)\bigr)=J'.
\end{equation}

Now for any $\mu$ satisfying \eqref{much} set
\begin{equation}
\label{tmu} T_\mu= J(L^2 - \mu),
\quad\dom\bigl(T_{\mu}\bigr)=\dom(L^2),
\end{equation}
and
\begin{equation}
\mathfrak{t}_\mu(x,y)=\lal LJx,Ly\ral-\mu\lal x, Jy\ral, \quad
x,y\in\dom(\mathfrak{t}_\mu)=\dom(L).
\end{equation}
Clearly, $\dom\bigl(T_{\mu}\bigr)\subset\dom(\mathfrak{t}_\mu)$ and
$\mathfrak{t}_\mu(x,y)=\lal x,T_{\mu} y\ral$ for any
$x,y\in\dom\bigl(T_{\mu}\bigr)$. Further, introduce the set
$\mathscr{E}$ of ordered orthonormal two-element systems in $\got H$
by
\begin{equation*}
\mathscr{E} = \bigl\{ \{e_{-},e_{+}\} \subset
\dom(\mathfrak{t}_\mu)\,\, \big| \quad e_{\pm}\in \got H_{\pm},
\quad \|e_{\pm}\| = 1\bigr\}.
\end{equation*}
Then by Lemma \ref{wt}
\begin{equation}
\label{inclusion}
\cW\bigl(T_{\mu}\bigr)\subset\bigcup_{\mathbf{e}\in \mathscr{E}}
\cW\bigl(\got t^{\mathbf{e}}_{\mu}\bigr),
\end{equation}
where  $\got t^{\mathbf{e}}_{\mu}$ are $2\times 2$ matrices given by
\begin{equation*}
%\label{Tbe}
\got t^{\mathbf{e}}_{\mu}=\left(\begin{array}{cc}
\mathfrak{t}_\mu(e_-,e_-) & \mathfrak{t}_\mu(e_-,e_+) \\
\mathfrak{t}_\mu(e_+,e_-) & \mathfrak{t}_\mu(e_+,e_+)
\end{array}\right),\quad{\mathbf e}=\{e_-,e_+\}\in\mathscr{E}.
\end{equation*}
By taking into account that $A\smile J$ and $V\frown J$, one
observes
\begin{equation}
\label{Tbee} \got t^{\mathbf{e}}_{\mu} =
\begin{pmatrix}
\mu - \| A e_{-}\|^2 - \|V e_{-}\|^2 &
-\overline{(\lal A e_{+}, V e_{-}\ral + \lal V e_{+},A e_{-}\ral)} \\
\lal A e_{+}, V e_{-}\ral + \lal V e_{+},A e_{-}\ral & \|A e_{+}\|^2
+ \|V e_{+}\|^2 - \mu
\end{pmatrix}.
\end{equation}

From \eqref{intervals1} it follows that for $\{e_{-},e_{+}\} \in
\mathscr{E}$
\begin{equation}
\label{Aepm} \| A e_{-} \| \le a \text{\, and \,} \| A e_{+}\| \ge
b.
\end{equation}
Hence, under the assumption \eqref{much} by Lemma \ref{ellipse} it
follows from \eqref{Tbee} and \eqref{Aepm} that  for all
$\mathbf{e}\in \mathscr{E}$ the numerical ranges $\cW\bigl(\got
t^{\mathbf{e}}_{\mu}\bigr)$ are elliptical discs that lie in the
open right half-plane $\{ z\in\bbC \, | \, \re z > 0 \}$. Then
\eqref{inclusion} implies that the numerical range
$\cW\bigl(T_\mu\bigr)$ also lies in the open right half-plane, that
is, the operator $T_{\mu}$ is strictly accretive. Now taking into
account \eqref{JLm} and \eqref{tmu}, Theorem \ref{inv} yields that
the involution $J$ and $J'$ are in the acute case. Moreover, for the
direct rotation $U$ from $J$ to $J'$ the following inequality holds
\begin{equation}
\label{sector2} \vartheta(U)\leq\frac{1}{2}\arctan k(T_{\mu}),
\end{equation}
where $\mu$ is an arbitrary point from the interval~\eqref{much}. In
its turn, inclusion \eqref{inclusion} implies
\begin{equation}
\label{supk} k(T_{\mu}) \leq \sup_{\mathbf{e}\in \mathscr{E}} \,
k(\got t^{\mathbf{e}}_{\mu}).
\end{equation}
Since
\begin{equation*}
|\lal A e_{+}, V e_{-}\ral + \lal V e_{+},A e_{-}\ral|\leq \|A
e_{+}\|\| V e_{-}\|+\|A e_{-}\|\|V e_{+}\|, \quad
\mathbf{e}=\{e_-,e_+\}\in\mathscr{E},
\end{equation*}
by  Lemma~\ref{ellipse} it follows  from \eqref{Tbee} that
\begin{equation}
\label{kt} k\bigl(\got t^\mathbf{e}_{\mu}\bigr) \le
f_\mu(\alpha_-,\alpha_+,v_-,v_+),
\end{equation}
where
\begin{equation*}
f_\mu(\alpha_-,\alpha_+,v_-,v_+)=\frac{\alpha_-v_+ + \alpha_+v_-}
{(\mu - \alpha_-^2-v_-^2)^{1/2} (\alpha_+^2 + v_+^2 - \mu)^{1/2}}
\end{equation*}
with $\alpha_\pm=\|Ae_\pm\|$ and $v_\pm=\|Ve_\pm\|$.

By \eqref{Aepm} we have
\begin{equation}
\label{alfpm} 0\leq\alpha_-\leq a\text{\, and \,}\alpha_+\geq b,
\end{equation}
while
\begin{equation}
\label{vepm} 0\leq v_-\leq\|V\|\text{\, and \,}0\leq v_+\leq\|V\|.
\end{equation}
A direct computation shows that the supremum of the function $f_\mu$
over the set in $\bbR^4$ constrained by \eqref{alfpm} and
\eqref{vepm} equals
\begin{equation*}
\varkappa(\mu) = \left\{ \begin{aligned} \frac{\|V\|(a+b)}{(\mu -
a^2 - \|V\|^2)^{1/2}(b^2+\|V\|^2 -\mu)^{1/2}}
\quad & \text{\, if \,} a(b^2 - \mu) > b \|V\|^2, \\
\frac{[b^2 \|V\|^2 + a^2 (b^2 - \mu)]^{1/2}} {(\mu - a^2 -
\|V\|^2)^{1/2}(b^2-\mu)^{1/2}} \quad &  \text{\, if \,} a(b^2 - \mu)
\le  b \|V\|^2.
\end{aligned}
\right.
\end{equation*}
Then by \eqref{sector2}--\eqref{kt} one infers that
\begin{equation*}
%\label{thekap}
\vartheta(U)\leq\frac{1}{2}\arctan \varkappa(\mu) \text{\, for any
\,} \mu\in(a^2 + \|V\|^2, b^2).
\end{equation*}
In particular,
\begin{equation}
\label{thekap1} \vartheta(U)\leq\frac{1}{2}\arctan
\varkappa_{\mathrm{min}},
\end{equation}
where
\begin{equation}
\label{kapmin} \varkappa_{\mathrm{min}}=\inf\limits_{a^2 +
\|V\|^2<\mu<b^2}\varkappa(\mu).
\end{equation}

By inspection, the function $\varkappa(\mu)$ is continuously
differentiable on the interval $(a^2 + \|V\|^2, b^2)$. The (global)
minimum of $\varkappa$ on this interval is just equal to
$\kappa\bigl(\|V\|\bigr)$. By \eqref{thekap1} the equality
$\varkappa_{\mathrm{min}}=\kappa\bigl(\|V\|\bigr)$ proves the bound
\eqref{thebest}.

The uniqueness of an involution $J'$ with the properties
$(i)$--$(iii)$ follows from Theorem \ref{Juniq1}. Estimate
\eqref{thprime2} is an immediate corollary to Theorem \ref{Juniq}.

The proof is complete.
\end{proof}

\begin{rem}
Notice that in the case where the operator $A$ is bounded, the
estimate $\|\sE_L(\sigma'_-)-\sE_A(\sigma_-)\|< \dfrac{\sqrt{2}}{2}$
(or equivalently $\vartheta(U)<{\pi}/{4}$) may be obtained by
combining \cite[Theorem 1 $(ii)$]{KMM3} and \cite[Theorem
5.6]{KMM5}.
\end{rem}

Theorem \ref{main0} is an immediate corollary to Theorem
\ref{trio1}.

\begin{proof}[Proof of Theorem \ref{main0}]
Let $\Delta$ again denote the finite gap of the set $\sigma_{+}$
that contains $\sigma_{-}$. Obviously, $|\Delta| \ge 2 d$ and thus
$\| V\| < d \le \sqrt{d (|\Delta| - d)}$. By Theorem \ref{trio1} one
concludes that
\begin{equation*}
\| \sE_{L}(\sigma'_{-}) - \sE_{A}(\sigma_{-})\| \le \sin \Bigl(
\frac{1}{2}\arctan \kappa\bigl(\|V\|\bigr) \Bigr)
\end{equation*}
with $\kappa(v)$ given by \eqref{kappaV}. Observing that for $0\leq
v<d$
\begin{equation*}
\kappa(v) \le \dfrac{2vd}{d^2 - v^2} = \tan\Bigl( 2 \arctan
\frac{v}{d} \Bigr)
\end{equation*}
completes the proof.
\end{proof}

\begin{ex}
\label{tsharp} Let $A$ be a self-adjoint operator on $\got H =
\mathbb{C}^4$ defined by
\begin{equation*}
A = \mathop{\mathrm{diag}}\{ -b, -a, a, b\}, \quad 0 \le a < b.
\end{equation*}
Divide the spectrum of $A$ into the two disjoint sets
$\sigma_-=\{-a,a\}$ and $\sigma_+=\{-b,b\}$. Clearly,
$d=\dist(\sigma_-,\sigma_+)=b-a>0$. The interval  $\Delta=(-b,b)$
appears to be the gap of the set $\sigma_+$ containing the set
$\sigma_-$. The involution $J = \sE_{A}(\sigma_{+}) -
\sE_A(\sigma_{-})$ reads
\begin{equation*}
J = \mathop{\mathrm{diag}}\{ +1, -1, -1, +1\}.
\end{equation*}
Assume that $V$ is a $4\times4$ matrix of the form
\begin{equation}
\label{Vform} V =
  \begin{pmatrix}
    0 & v_1 & v_2 & 0 \\
    v_1 & 0 & 0 & v_2 \\
    v_2 & 0 & 0 & v_1 \\
    0   & v_2 & v_1 & 0
  \end{pmatrix},
\end{equation}
where $v_1,v_2 \ge 0$. By inspection, $V$ anticommutes with $J$ and
$\|V\| = v_1 + v_2$. The involution $J'=\sE_L(\mathbb{R}\backslash
\Delta) - \sE_L(\Delta)$ is computed explicitly as soon as the
eigenvectors of the $4\times4$ matrix {$L=A+V$} are found. Under the
assumption that \eqref{Vdd} holds, that is, for $\|V\|^2 < b^2 -
a^2$, the explicit evaluation of the spectral angle of the direct
rotation $U$ from $J$ to $J'$ results in
\begin{equation*}
\vartheta(U) = \dfrac{1}{2}\arctan\left(\frac{2 a(v_1 - v_2) + 2
b\|V\|} {b^2 - a^2 - \|V\|^2 + (v_1-v_2)^2}\right).
\end{equation*}
Taking into account that the value of $v_1-v_2$ for different
matrices \eqref{Vform} with the same norm $\|V\|$ runs through the
interval $[-\|V\|,\|V\|]$, one easily verifies that the maximal
possible value $\vartheta_{\rm max}$ of $\vartheta(U)$ is equal to
\begin{equation}
\label{tetamax} \vartheta_{\rm max}=\frac{1}{2}\arctan\kappa(\|V\|),
\end{equation}
with $\kappa(v)$ given by \eqref{kappaV}. In particular, if $a=0$
then
\begin{equation}
\label{tetamax0} \vartheta_{\rm
max}=\arctan\left(\frac{\|V\|}{d}\right).
\end{equation}

\end{ex}

\begin{rem}
Example \ref{tsharp} shows the following.
\begin{itemize}
\item[$(i)$] Estimate
\eqref{thebest} of Theorem \ref{trio1} is sharp. This is proven by
equality \eqref{tetamax}.
\item[$(ii)$] Estimate \eqref{edif1} of Theorem \ref{main0} is also sharp.
This is proven by equality \eqref{tetamax0}.
\end{itemize}
\end{rem}
\begin{rem}
\label{conj} We \textit{conjecture} that estimate \eqref{edif1} of
Theorem \ref{main0} also holds for $d\leq \|V\|<\sqrt{2}d$.
\end{rem}

% ------------------------------------------------------------------------
\subsection*{Acknowledgment}
This work was supported by the Deutsche Forschungsgemeinschaft
(DFG), the Heisenberg-Landau Program, and the Russian Foundation for
Basic Research.
% ------------------------------------------------------------------------


\begin{thebibliography}{99}

\bibitem{AdLMSr} V.\,M. Adamyan,  H. Langer, R. Mennicken, and
J. Saurer, \textit{Spectral components of selfadjoint block operator
matrices with unbounded entries}, Math. Nachr. \textbf{178} (1996),
43 -- 80.

\bibitem{ALT} V. Adamyan, H. Langer, and C. Tretter,
\textit{Existence and uniqueness of contractive solutions of some
Riccati equations}, J. Funct. Anal. \textbf{179} (2001) 448 -- 473.

\bibitem{AMM} S. Albeverio, K. A. Makarov and A. K. Motovilov,
\textit{Graph subspaces and the spectral shift function}, Canadian
Journal of Mathematics, {\bf 55}:3 (2003),  449 -- 503;
arXiv: math.SP/0105142.

\bibitem{AI} T. Y. Azizov and I. S. Iokhvidov,
\textit{Linear Operators in Spaces with an Indefinite Metric}, John
Wiley \& Sons, Chichester, 1989.

\bibitem{Birman-Solomjak} M. Sh. Birman and M. Z. Solomjak,
\textit{Spectral Theory of Self-Adjoint Operators on a Hilbert
Space}, Leningrad State University, Leningrad, 1980 (Russian).

\bibitem{Davis} C. Davis, \textit{Separation of two linear subspaces},
Acta Sci. Math. Szeged \textbf{19} (1958), 172 -- 187.

\bibitem{Davis:123} C.~Davis, \textit{The rotation of
eigenvectors by a perturbation.  I and II}, J. Math. Anal. Appl.
\textbf{6} (1963), 159 -- 173; \textbf{11} (1965), 20 -- 27.

\bibitem{DK} C.~Davis and W.~M.~Kahan, \textit{The
rotation of eigenvectors by a perturbation. III}, SIAM J. Numer.
Anal. \textbf{7} (1970), 1 -- 46.

\bibitem{Rao} K. E. Gustafson and D. K. M. Rao, \textit{Numerical Range},
Springer, N. Y., 1997.

\bibitem{Halmos:69} P. R. Halmos, \textit{Two subspaces}, Trans.
Amer. Math. Soc. \textbf{144} (1969), 381--389.

\bibitem{HMM} V. Hardt, R. Mennicken, and A. K. Motovilov,
\textit{Factorization theorem for the transfer function associated
with a $2\times2$ operator matrix having unbounded couplings}, J.
Oper. Th. \textbf{\bf 48:1} (2002), 187 -- 226.

\bibitem{Kato} T.~Kato, \textit{Perturbation Theory for Linear
Operators}, Springer--Verlag, Berlin, 1966.

\bibitem{KMM1} V. Kostrykin, K. A. Makarov,
and A. K. Motovilov, \textit{On a subspace perturbation problem},
Proc. Amer. Math. Soc. {\bf 131} (2003), 3469 -- 3476;
arXiv: math.SP/0203240.

\bibitem{KMM2} V. Kostrykin, K. A. Makarov,
and A. K. Motovilov, \textit{Existence and uniqueness of solutions
to the operator Riccati equation. A geometric approach},
Contemporary Mathematics (AMS) {\bf 327} (2003), 181 -- 198;
arXiv: math.SP/0207125.

\bibitem{KMM3} V. Kostrykin, K. A. Makarov,
and A. K. Motovilov, \textit{On the existence of solutions to the
operator Riccati equation and the tan\,$\Theta$ theorem},  Int. Eq.
Op. Th. {\bf 51} (2005), 121 -- 140;
arXiv:\,math.SP/0210032 v2.

\bibitem{KMM4} V. Kostrykin, K. A. Makarov,
and A. K. Motovilov, \textit{Perturbation of spectra and spectral
subspaces}, Trans. Amer. Math. Soc. (to appear);
arXiv:\,math.SP/0306025 v1.

\bibitem{KMM5} V. Kostrykin, K. A. Makarov,
and A. K. Motovilov, \textit{A generalization of the $\tan 2\Theta$
Theorem}, Operator Theory: Adv. Appl.  \textbf{149} (2004), 349 --
372; arXiv:\,math.SP/0302020.

\bibitem{MenShk} R. Mennicken and A. A. Shkalikov,
\textit{Spectral decomposition of symmetric operator matrices},
Math. Nachr. \textbf{179} (1996), 259 -- 273.

\bibitem{SzNagy} F. Riesz and B. Sz.-Nagy, \textit{Le\c{c}ons d'analyse
fonctionelle}, 2nd ed., Acad\'emiai Kiado, Budapest, 1953.

\end{thebibliography}
\end{document}